\documentclass[10pt]{article}
\usepackage{amssymb,amsfonts,amsmath,amsthm,amsbsy}%
\usepackage[dvips]{graphicx,color}%
\theoremstyle{plain}%
\newtheorem{theorem}{Theorem}[section]%
\newtheorem{observation}[theorem]{Observation}%

\newtheorem{lemma}[theorem]{Lemma}%
\newtheorem{conjecture}[theorem]{Conjecture}%
\newtheorem{remark}[theorem]{Remark}%
\newenvironment{Proof}[1][.]%
 {\begin{trivlist}\item[]\textbf{Proof#1 }}%
 {\hspace*{\fill}$\rule{0.3\baselineskip}{0.35\baselineskip}$\end{trivlist}}
\newcommand{\fpu}{\mathrm{FPU}}

\newcommand{\R}{\mathbb{R}}

\newcommand{\calD}{\mathcal{D}}
\newcommand{\calJ}{\mathcal{J}}

\newcommand{\calH}{\mathcal{H}}

\newcommand{\calR}{\mathcal{R}}
\newcommand{\calS}{\mathcal{S}}
\newcommand{\calW}{\mathcal{W}}
\newcommand{\calO}{\mathcal{O}}

\newcommand{\bal}{{\overline \alpha}}
\newcommand{\bt}{{\overline  t}}

\newcommand{\rmb}{\mathrm{b}}

\newcommand{\rmd}{\mathrm{d}}

\newcommand{\rmf}{\mathrm{f}}

\newcommand{\rmrh}{\mathrm{rh}}

\newcommand{\sgn}{\mathrm{sgn}}

\newcommand{\eps}{{\varepsilon}}

\newcommand{\lst}{{\mathrm{L}}}
\newcommand{\mst}{{\mathrm{M}}}
\newcommand{\rst}{{\mathrm{R}}}
\newcommand{\bigpar}{\par\quad\par}
\newcommand{\micro}{{\rm \,FPU}}
\newcommand{\cons}{{\rm \,cons}}

\newcommand{\SpTi}{space-time}%
\newcommand{\tdots}{{...}}%

\newcommand{\pair}[2]{{\left({#1},\,{#2}\right)}}
\newcommand{\bpair}[2]{{\big({#1},\,{#2}\big)}}
\newcommand{\at}[1]{{\left({#1}\right)}}
\newcommand{\nat}[1]{{({#1})}}
\newcommand{\bat}[1]{{\big(#1\big)}}
\newcommand{\Bat}[1]{{\Big(#1\Big)}}
\newcommand{\triple}[3]{{\left({#1},\,{#2},\,{#3}\right)}}

\newcommand{\abs}[1]{{\left|{#1}\right|}}

\newcommand{\MaFld}[1]{\overline{#1}}

\newcommand{\MiTime}{t}
\newcommand{\MiSpace}{x}

\newcommand{\MiLagr}{\alpha}

\newcommand{\MaSpace}{\MaFld{\MiSpace}}

\newcommand{\MaTime}{\MaFld{\MiTime}}
\newcommand{\MaLagr}{\MaFld{\MiLagr}}

\newcommand{\MaLagrDer}[1]{\partial_{\,\MaLagr\,}{#1}}

\newcommand{\MaTimeDer}[1]{\partial_{\,\MaTime\,}{#1}}
\newcommand{\dint}[1]{\,\mathrm{d}#1}
\newcommand{\jump}[1]{{|\![#1]\!|}}

\newcommand{\meas}[1]{{\left\langle #1\right\rangle }}

\newcounter{Lcount}
\def\enum{\setcounter{Lcount}{1}
\begin{list}{\arabic{Lcount}.}{\usecounter{Lcount}\leftmargin=0.5em}}
\newcommand{\figdraft}{false}%
\newcommand{\figdirectory}{pictures}%
\newcommand{\figfile}[1]{\figdirectory/#1}%
\newlength{\mhpicDwidth}
\newlength{\mhpicDvsep}
\newlength{\mhpicPwidth}
\newlength{\mhpicPhsep}
\usepackage{float}%
\numberwithin{theorem}{section}%
\numberwithin{equation}{section}%
\numberwithin{figure}{section}
\begin{document}%
%
%
\title{
Riemann solvers and undercompressive shocks\\ of convex FPU chains
}%
\date{\today}%
\author{Michael Herrmann\thanks{ %
    University of Oxford,
    Mathematical Institute, Centre for Nonlinear PDE (OxPDE),
    24-29 St Giles', Oxford OX1 3LB, England, michael.herrmann@maths.ox.ac.uk.}
\and
Jens D.M. Rademacher\thanks{%
    National Research Centre for Mathematics and Computer Science (CWI, MAS),
    Science Park 123, 1098 XG Amsterdam, the Netherlands,  rademach@cwi.nl.}
}
\maketitle
%
%
%
\begin{abstract}%
We consider FPU-type atomic chains with general convex potentials.
The naive continuum limit in the hyperbolic space-time scaling is
the p-system of mass and momentum conservation. We systematically
compare Riemann solutions to the p-system with numerical solutions
to discrete Riemann problems in FPU chains, and argue that the
latter can be described by modified p-system Riemann solvers. We
allow the flux to have a turning point, and observe a third type of
elementary wave (conservative shocks) in the atomistic
simulations. These waves are heteroclinic travelling waves and
correspond to non-classical, undercompressive shocks of the
p-system. We analyse such shocks for fluxes with one or more turning
points.

Depending on the convexity properties of the flux we propose FPU-Riemann 
solvers. Our numerical simulations confirm that Lax-shocks
are replaced by so called dispersive shocks. For convex-concave flux
we provide numerical evidence that convex FPU chains follow the
p-system in generating conservative shocks that are supersonic. For
concave-convex flux, however, the conservative shocks of the
p-system are subsonic and do not appear in FPU-Riemann solutions.
\end{abstract}%

%
%
%
%
%
%
%
\section{Introduction}\label{chIntro}
%
The derivation of effective continuum descriptions for
high-dimensional discrete systems is a fundamental tool for model
reduction in the sciences. Hamiltonian lattices, such as atomic
chains, naturally lead to nonlinear systems of conservation laws
which describe the leading order dynamics on the hyperbolic
space--time scale.  It is customary to neglect the higher
order terms and to study the leading order system by itself. The
rigorous mathematical validity of this reduction is a notoriously
difficult task and there are surprisingly few successes reported in
the literature.
\par
\begin{figure}[ht!]%
  \centering{%
  \includegraphics[width=0.7\textwidth, draft=\figdraft]%
  {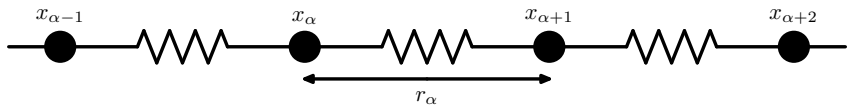}%
  }%
  \caption{The atomic chain with nearest neighbour interaction.}%
  \label{IntroFig1}%
\end{figure}

This paper concerns the macroscopic description of monoatomic chains
with nearest neighbour interactions, see Figure~\ref{IntroFig1}. In
their seminal paper \cite{FPU55} Fermi, Pasta and Ulam chains studied
such chains for interaction potential $\Phi$ whose non-harmonic part
involves only cubic or quartic terms. We allow for general convex
interaction potentials, but still refer to the systems as FPU chains.
\par
An important building block for the macroscopic descriptions of
FPU chains are the solutions to Riemann initial data in the
hyperbolic continuum limit. In this limit space and time are scaled
in the same way, and the amplitude of solutions is unconstrained.
Since the numerical study of Holian and Straub \cite{HS78} and from
rigorous results for the integrable Toda chain it is known that the
solutions to atomistic Riemann problems with either convex or
concave flux $\Phi^\prime$ obey a self-similar structure on the
hyperbolic scale: each solution consists of at most two elementary
waves that are separated by constant states. Each of these
elementary wave is either a rarefaction wave, in which the atomic
data vary smoothly on the macroscopic scale, or a \emph{dispersive
shocks}, in which strong microscopic oscillations spread out in
space and time.
\par
The starting point for our investigation was the observation that
for certain $\Phi$ a third kind of elementary waves can be observed
in FPU-Riemann problems. These waves, which we refer to as
conservative shocks, involve no oscillations and look like
`shocks' (jump discontinuities) on the macroscopic scale. Of course,
these waves are not exact shocks as they exhibit a
transition layer on the atomistic scale, but this layer is very
small and disappears in the hyperbolic scaling. To our knowledge the
appearance of conservative shocks in FPU-Riemann problems
was never reported before.
\par
\begin{figure}[ht!]
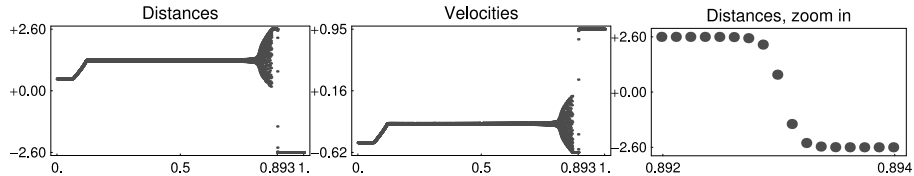

\centering{%
\setlength{\tabcolsep}{0cm}%
\begin{tabular}{ccc}%
\includegraphics[width=\mhpicDwidth, draft=\figdraft]%
{\figfile{rp_supersonic/rp_super_dist_mi}}%
&
\includegraphics[width=\mhpicDwidth, draft=\figdraft]%
{\figfile{rp_supersonic/rp_super_vel_mi}}%
&
\includegraphics[width=\mhpicDwidth, draft=\figdraft]%
{\figfile{rp_supersonic/rp_super_dist_zoom}}%
\end{tabular}%
}%
\caption{%
Riemann problems in convex FPU chains can involve three kinds of
elementary waves: rarefaction waves, dispersive shocks, and
conservative shocks. The first two pictures show snapshots of the
atomic distances and velocities against the scaled particle index
$\alpha/N$ for a chain with $N=8000$ particles and flux
$\Phi^\prime\at{r}=\tfrac{1}{16}\at{r+\arctan{r}}$; the third
pictures magnifies the region around the conservative
shock.}\label{f:intro_rp_supersonic}
\end{figure}

An illustrative example of an FPU-Riemann problem is plotted in
Figure~\ref{f:intro_rp_supersonic} and involves all types of
elementary waves (from left to right): a rarefaction wave, a
dispersive shock, and conservative shock.
\par
The naive approach to describe the FPU dynamics on the hyperbolic
scale assumes long-wave-length motion without microscopic
oscillations. Under this assumption one readily derives
the `p-system', which consists of the conservation laws for mass and
momentum in Lagrangian coordinates.
\par
It is well known, that the naive continuum limits of nonlinear
dispersive lattices provide a reasonable macroscopic model as long
as the macroscopic fields are smooth, see, e.g.,
\cite{Lax86,GL88,HL91,HLM94,LL96}. The nonlinearity, however,
usually causes shock phenomena, and the naive continuum limit fails
in this case. Instead, lattice systems like FPU typically produce
\emph{dispersive shocks}, in which the atoms self-organize into
strong microscopic oscillations. On the macroscopic scale such
dispersive shocks can be regarded as measure-valued solutions to the
naive continuum limit, see \eqref{e:YMMacroConsLaws} in
\S\ref{s:Found}.
\par
\begin{figure}[ht!]
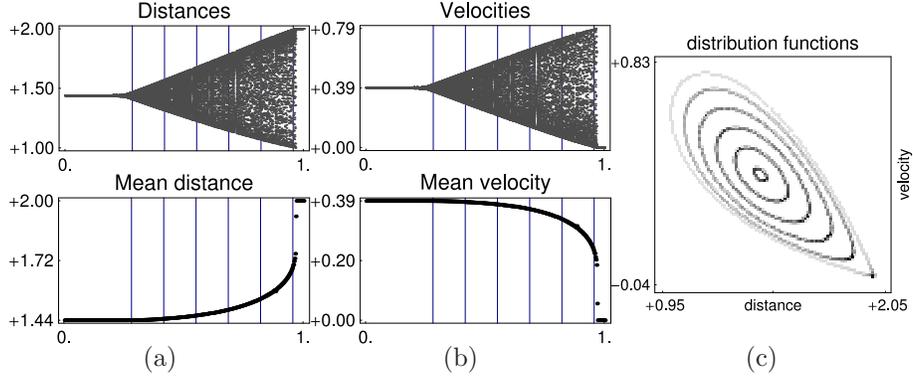
%
\centering{%
\setlength{\tabcolsep}{0cm}%
\begin{tabular}{ccc}
\begin{minipage}[c]{\mhpicDwidth}%
\includegraphics[width=\mhpicDwidth, draft=\figdraft]%
{\figfile{structure_dispersive_shocks/sds_dist_mi.eps}}%
\\%
\includegraphics[width=\mhpicDwidth, draft=\figdraft]%
{\figfile{structure_dispersive_shocks/sds_dist_ma.eps}}%
\end{minipage}%
&%
\begin{minipage}[c]{\mhpicDwidth}%
\includegraphics[width=\mhpicDwidth, draft=\figdraft]%
{\figfile{structure_dispersive_shocks/sds_vel_mi.eps}}%
\\%
\includegraphics[width=\mhpicDwidth, draft=\figdraft]%
{\figfile{structure_dispersive_shocks/sds_vel_ma.eps}}%
\end{minipage}%
&%
\begin{minipage}[c]{\mhpicDwidth}
\includegraphics[width=\mhpicDwidth, draft=\figdraft]%
{\figfile{structure_dispersive_shocks/sds_df_dist_vel.eps}}%
\end{minipage}%
\\(a) & (b) & (c)%
\end{tabular}%
}%
\caption{%
  Dispersive shocks arise naturally in FPU-Riemann problems
  as a consequence of energy conservation. Each dispersive shock is build of a
  one-parameter family of wave trains with
  a single soliton at the leading front. (a,b) Snapshots of atomic
  distances and velocities, and their local mean values.  (c) Superposition of
  several local
  distribution functions within the shock; positions of the mesoscopic
  {\SpTi} windows are marked by vertical lines in (a,b).
%
}
\label{f:shock}%
\end{figure}

The formation of dispersive shocks is a characteristic property of
Hamiltonian `zero dispersion' limits and a direct consequence of the
conservation of energy, see \S\ref{s:Found}. Moreover, it is known
from numerical studies and rigorous results for integrable systems,
see \cite{GP73,LL83,V85,Lax86,Lax91,LLV93,Kam00}, that the
oscillations in a dispersive shocks are modulated wave trains (period
travelling waves). Figure~\ref{f:shock} presents a typical example of
a dispersive shock in FPU chains. At the shock front, where the
amplitudes of the oscillations become maximal, the wave trains
converge to a supersonic soliton, that is a homoclinic travelling
wave.
\bigpar By combining non-classical hyperbolic theory of the p-system
with macroscopic theory, travelling waves and numerical observations
of FPU chains we characterise FPU-Riemann solvers for
oscillation-free initial data and fluxes with one turning point.
%
%
\subparagraph{Conservative shocks in FPU chains and the p-system.}
%
%
As shown in Figure~\ref{f:intro_rp_supersonic}, the solution to
FPU-Riemann problems with convex $\Phi$ can involve conservative
shocks. Below in \S\ref{s:Found} we explain that these waves
correspond to certain shocks in the p-system, namely those that
conserve the energy exactly. Among all p-system shocks, the set of
conservative shocks is quite small and if $\Phi^\prime$ has no
turning point conservative shocks cannot occur at all. The
conservative shocks in the p-system are non-classical shocks as they
violate the Lax-condition: for convex-concave $\Phi^\prime$ there
are \emph{fast undercompressive}, and hence \emph{supersonic} with
respect to both the left and right state; for concave-convex
$\Phi^\prime$ the conservative shocks are \emph{slow
undercompressive} and hence \emph{subsonic}.
\par
\begin{figure}[ht!]
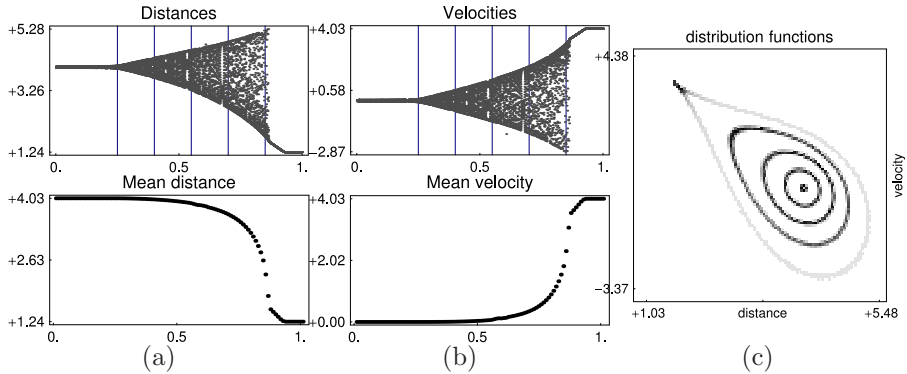

\centering{%
\setlength{\tabcolsep}{0cm}%
\begin{tabular}{ccc}%
\begin{minipage}[c]{\mhpicDwidth}%
\includegraphics[width=\mhpicDwidth, draft=\figdraft]%
{\figfile{qp2_composite_wave/qp2_cw_dist_mi.eps}}%
\\%
\includegraphics[width=\mhpicDwidth, draft=\figdraft]%
{\figfile{qp2_composite_wave/qp2_cw_dist_ma.eps}}%
\end{minipage}%
&%
\begin{minipage}[c]{\mhpicDwidth}%
\includegraphics[width=\mhpicDwidth, draft=\figdraft]%
{\figfile{qp2_composite_wave/qp2_cw_vel_mi.eps}}%
\\%
\includegraphics[width=\mhpicDwidth, draft=\figdraft]%
{\figfile{qp2_composite_wave/qp2_cw_vel_ma.eps}}%
\end{minipage}%
&%
\hspace{\mhpicPhsep}%
\begin{minipage}[c]{\mhpicPwidth}%
\includegraphics[width=\mhpicPwidth, draft=\figdraft]%
{\figfile{qp2_composite_wave/qp2_cw_df_dist_vel.eps}}%
\end{minipage}%
\\(a)&(b)&(c)%
\end{tabular}%
}%
\caption{%
FPU-Riemann problem with initial data corresponding to a subsonic
conservative shock in the p-system: Instead of a conservative
shock the atoms self-organize into a dispersive shock with attached
rarefaction wave. The leading soliton in the dispersive shock is
hence {sonic}. (a) atomic distances;  (b) atomic velocities; (c)
five local distributions functions corresponding to the vertical
lines in the snapshots.
%
}\label{f:consShock.B}
\end{figure}

We numerically discovered that conservative shocks occur naturally
in FPU-Riemann problems chains if they are supersonic. However,
numerical simulations with concave-convex $\Phi^\prime$ never
generated anything close to a jump discontinuity. Instead, we
typically observe solutions as plotted in
Figure~\ref{f:consShock.B}: the solution appears to be a composite
wave of a dispersive shock with \emph{attached} rarefaction wave.
\par
Conservative shocks in FPU chains are naturally related to atomistic
fronts which are `heteroclinic' travelling waves. A bifurcation result
of Iooss \cite{Ioo00} shows the existence of small amplitude atomistic
fronts for convex-concave flux (in which case the conservative shocks
are supersonic). The authors have improved this result in \cite{HR09}
by showing that each front must correspond to a conservative shock in
the p-system and that supersonic shocks with arbitrary large jump
height can be realised by an atomistic front.
%
%
%
\subparagraph{FPU-Riemann solvers.}
%
A main goal of this paper is the derivation of FPU-Riemann solvers
which predict the number and the type of the elementary waves that
result from arbitrary Riemann initial data. To this end we
systematically compare numerical simulation for various FPU chains
with certain Riemann solvers for the p-system. The FPU-Riemann
solver for the Toda chain is well understood, see for instance
\cite{Kam93,DM98}, but there is no complete picture for
non-integrable chains as the available numerical studies only
concern special types of Riemann initial data, such as the `piston
problem' in \cite{HS78}. Moreover, we are not aware of any previous
analytical or numerical investigation of any FPU-Riemann problem
which allows for conservative shocks.
\bigpar
In the classical case with either convex or concave flux
$\Phi^\prime$, the numerical simulations indicate that the solution
to each FPU-Riemann problem can be described by an adapted classical
solver, in which Lax-shocks are replaced by dispersive shocks. More
precisely, from each given left state there emanate four curves
(wave sets) in the state space. Two of them correspond to
rarefaction waves and appear also in each solver for the p-system.
Instead of Lax-shock curves, however, we find two
dispersive shock curves. The solution to each FPU-Riemann problem is
then completely determined by these wave sets. In particular, in the
classical case each Riemann solution consists (generically) of a
left moving $1$-wave and a right moving $2$-wave, where each wave is
either a rarefaction wave or a dispersive shock.
\par%
In presence of turning points of $\Phi^\prime$ the Riemann solvers for
both FPU chains and the p-system are more complicated because now the
solutions to general Riemann problems involve \emph{composite waves},
which consist of two elementary waves from the same family, and may
also involve undercompressive shocks. The Riemann solvers in these
cases can be described in terms of modified wave sets, but for the
p-system they are not unique. Consequently, different p-system Riemann
solvers are possible, such as the `conservative' and the `dissipative'
solver described in \S\ref{s:riem}. For FPU chains, however, the
underlying atomistic dynamics determines the solver uniquely for each
$\Phi$.
\bigpar
\begin{figure}[ht!]
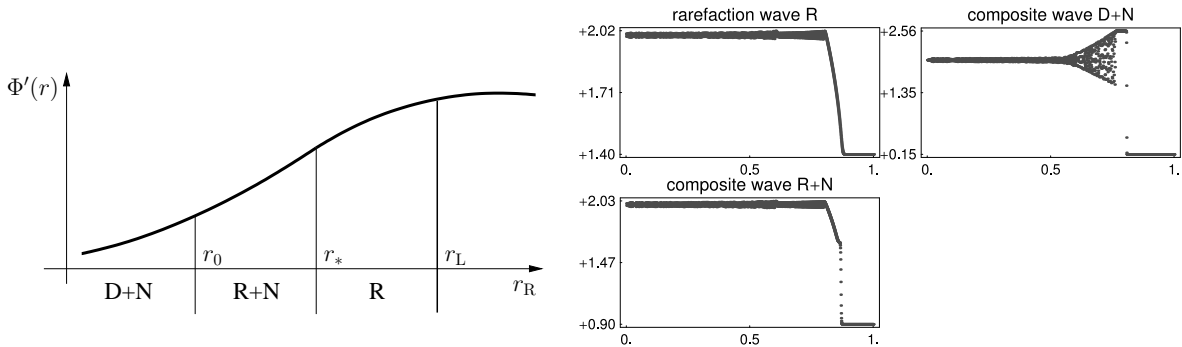
%
\begin{center}%
\begin{tabular}{cc}
\begin{minipage}[c]{0.45\textwidth}
\includegraphics[width=\textwidth, draft=\figdraft]%
{\figfile{xfig/riemann_2/fpu_mod_rare_curve}}%
\end{minipage}
&
\begin{minipage}[c]{2\mhpicDwidth}%
\begin{minipage}[t]{\mhpicDwidth}%
\includegraphics[width=\mhpicDwidth, draft=\figdraft]%
{\figfile{fpu_solver/fpusol_r.eps}}\\
\includegraphics[width=\mhpicDwidth, draft=\figdraft]%
{\figfile{fpu_solver/fpusol_rn.eps}}\\
\end{minipage}%
\begin{minipage}[t]{\mhpicDwidth}%
\includegraphics[width=\mhpicDwidth, draft=\figdraft]%
{\figfile{fpu_solver/fpusol_dn.eps}}\\
\end{minipage}%
\end{minipage}%
\end{tabular}
\end{center} %
\caption{Sketch of the modified rarefaction curve
($r_\rst\leq{r_\lst}$) for convex-concave $\Phi^\prime$ and given
$r_\lst$, and corresponding FPU-Riemann solutions: For
$r_\rst\lesssim{r}_\lst$ one still finds a rarefaction wave (R), but
if $r_\lst$ crosses the turning point $r_*$ of $\Phi^\prime$ a
non-classical, a conservative shock (N) nucleates. If $r_\rst$
decreases further the rarefaction wave is replaced by a dispersive
shock (D).
}%
\label{f:fpu_solver}
\end{figure}

In \S\ref{s:riem} we investigate systematically numerical solutions to
FPU-Riemann problems with convex-concave and concave-convex
$\Phi^\prime$ and argue that the solutions can be described by adapted
p-system solvers. More precisely, in the convex-concave we propose to
adapt the conservative solver, which predicts composite waves with
supersonic conservative shocks. The modified rarefaction curve of such
a solver is illustrated in Figure~\ref{f:fpu_solver}.  In the
concave-convex case, however, the simulations indicate that FPU-Riemann
solutions can be described by an adapted dissipative solver, whose
modified wave sets do not involve conservative shocks.
%
%
In both non-classical FPU solvers the Lax shocks are replaced by
dispersive shocks, and in the nucleation criterion for composite
waves the dispersive shock front velocity plays the role of the
Rankine-Hugeniot condition. These adaption rules lead to an adequate
description of the numerical solution to FPU-Riemann problems.
Moreover, the difference between Rankine-Hugeniot and dispersive
shock front velocity provides an explanation for the absence of
conservation shock in the concave-convex case: the nucleation
criterion for conservative shocks cannot be satisfied. However, it is
not known how to predict the shock front velocity from the initial data.
\bigpar
We believe that these results give new structural insight into the
hyperbolic nature of the continuum limit and the role of the p-system
in it. Moreover, they open up new avenues for analytical
investigation, some of which we phrase in the form of conjectures.
\par
We emphasize that the situation dramatically changes for
fluxes with more than one turning point. On the p-system level,
non-classical Riemann solvers can still be constructed, but these no
longer provide a basis for the continuum limit of the FPU chain. The
reason is that we numerically find energy conserving shocks between
wave trains. This is a new phenomenon  but its study is beyond the
scope of this article. Another problem left for future investigation
concerns Riemann problems where oscillations in form of wave trains
are already imposed in the initial data. One of the expected
new phenomena in this case is the onset of two-phase oscillations.
Finally, it would be interesting to study 
cold initial data with more than one jump discontinuity, 
because then two-phase oscillations can be created by the interaction of two 
dispersive shocks.
%
%
%
\subparagraph{Organisation of the paper.}
In \S\ref{s:Found} we collect some facts about FPU chains and the
p-system; we especially discuss how oscillatory FPU data give rise
to measure-valued solutions to the p-systems and briefly outline the
concept of  modulated wave trains. \S\ref{s:num} concerns the
numerical simulation of FPU-Riemann problems and contains our
observations about self-similarity on the macroscopic scale; we also
investigate the fine-structure of dispersive and conservative
shocks. \S\ref{s:riem} is devoted to FPU-Riemann problems. We start
with numerical results for the classical case with either convex or
concave $\Phi^\prime$. Afterwards we discuss the conservative and
dissipative solvers for the p-system, and proceed with FPU-Riemann
problems in the non-classical cases. Finally, in
\S\ref{s:mathconsshock} we prove some analytical results about
conservative shocks in the p-system.
%
%
%
%
%
\section{Preliminaries of FPU chains and the p-system}\label{s:Found}
%
%
%
Convex FPU chains consist of $N$ identical particles which are
nearest neighbour coupled in a convex potential $\Phi:\R\to \R$ by
Newton's equations
\begin{equation}
\label{e:NewEqn1}%
\ddot x_\alpha = \Phi'(x_{\alpha+1} - x_\alpha) - \Phi'(x_\alpha -
x_{\alpha-1}),
\end{equation}
where $\dot{~} = \frac{\rmd}{\rmd t}$ is the time derivative,
$x_\alpha(t)$ the atomic position, $\alpha=1,\ldots,N$ the particle
index. We consider nonlinear force $\Phi'$, referred to as the flux.
A prominent example for nonlinear $\Phi^\prime$ (without turning
points) is the completely integrable Toda chain, see
\cite{Tod70,Fla74,Hen74}, with
\begin{equation}
\label{e:TodaPot}%
\Phi_\mathrm{Toda}(r)=\exp\at{1-r}-(1-r).
\end{equation}
\par
For our purposes it is convenient to use the atomic distances
$r_\MiLagr=\frac{x_{\MiLagr+1}-x_\MiLagr}{(\alpha+1)-\alpha}$ and
velocities $v_\MiLagr=\dot{x}_\MiLagr$ as the basic variables,
changing (\ref{e:NewEqn1}) to the system
\begin{align}%
\label{e:NewEqn2} \dot{r}_\MiLagr = v_{\MiLagr+1}-v_\MiLagr\;,\quad
\dot{v}_\MiLagr = \Phi^{\prime}\at{r_\MiLagr} -
\Phi^{\prime}\at{r_{\MiLagr-1}}.
\end{align}
Note that while the Toda potential, and also the potentials used
below, allows for negative distances, this is essentially a matter
of suitably shifting the minimum by $\Phi \to \Phi(\cdot + r_0)$.
\par
We are interested in the thermodynamic limit $\varepsilon=1/N\to 0$
in the \emph{hyperbolic} scaling of the \emph{microscopic}
coordinates $\MiTime$ and $\MiLagr$. This scaling is defined by the
\emph{macroscopic} time $\MaTime = \varepsilon\MiTime$ and particle
index $\MaLagr = \varepsilon\MiLagr$. It is natural to scale the
atomic positions in the same way, i.e.
$\MaSpace=\varepsilon\MiSpace$, which leaves atomic distances and
velocities scale invariant.
%
%
In the limit $\varepsilon=0$ the spatial variable $\MaLagr$ becomes
continuous and the high dimensional ODE \eqref{e:NewEqn1} should be
replaced by a continuum limit, i.e., by a system of a few
macroscopic PDEs. Microscopic oscillations can be naturally
interpreted as a form of \emph{temperature} in the chain, see
\cite{AMSMSP:DHR}, and accordingly we refer to oscillation-free
limits as \emph{cold}.
%
%
%
\subparagraph{Evolution of cold data and the p-system}
%
To derive the p-system as the simplest model for the macroscopic
evolution of FPU chains we assume macroscopic fields
$r(\MaTime,\MaLagr)$ and $v(\MaTime,\MaLagr)$ such that
$r_\MiLagr\at\MiTime=
r\pair{\varepsilon\MiTime}{\varepsilon\MiLagr}$,
$v_\MiLagr\at\MiTime=%
v\pair{\varepsilon\MiTime}{\varepsilon\MiLagr}$.  This ansatz
corresponds to \emph{cold motion} as it assumes that there are no
microscopic oscillations in the chain.  Substitution into
(\ref{e:NewEqn2}) and taking the limit $\varepsilon\rightarrow0$
yields the macroscopic conservation laws for mass and momentum
\begin{align}%
\label{e:ColdModEqn} %
\MaTimeDer{r}-%
\MaLagrDer{v}=0,\quad
\MaTimeDer{v}-%
\MaLagrDer{\Phi^\prime\at{r}}%
&=0.%
\end{align}
It is well known that the p-system is hyperbolic for convex $\Phi$ and
that for smooth solutions the energy is conserved via
\begin{align}
\notag
\MaTimeDer{}\at{\textstyle{\frac{1}{2}v^2}+\Phi\at{r}}-%
\MaLagrDer{}\at{v\,\Phi^\prime\at{r}}=0.%
\end{align}
In the p-system a shock propagates with a constant shock speed
$c_\rmrh$ so that $r$ and $v$ satisfy the Rankine-Hugeniot jump
conditions for mass and momentum
\begin{align}%
\label{e:ColdDataJump}%
c_\rmrh\jump{r}+\jump{v}=0,\quad
c_\rmrh\jump{v}+\jump{\Phi^\prime\at{r}}=0,
\end{align}
where $\jump{x} = x_\lst-x_\rst$ denotes the jump. The main
observation is that for either convex or concave flux $\Phi'$ the
jump conditions \eqref{e:ColdDataJump} imply that the jump condition
for the energy must be violated, i.e.,
\begin{align}%
\label{e:ColdDataEnJump}
c_\rmrh\jump{\tfrac{1}{2}v^2+\Phi\at{r}}+%
\jump{v\,\Phi^\prime\at{r}}\neq0,
\end{align}
see Theorem~\ref{t:consData} below. Here the Lax criterion selects
the shocks with negative production.
%

\subparagraph{Onset of dispersive shocks in FPU chains}
%
%
It is known that the p-system provides a reasonable thermodynamic
limit for FPU chains in the following sense:
Preparing cold FPU initial data with smooth profile functions $r$
and $v$, the atomistic dynamics reproduces a solution to the
p-system provided that the latter has a smooth solution. 
This can be understood as a manifestation of Strang's Theorem \cite{Str64};
we refer to \cite{DH08} for numerical simulations and 
to \cite{GL88,HL91} for a similar discussion in the context of
other lattice systems.
\par
At some
critical time, however, the p-system forms a shock; this shock still
conserves mass and momentum according to \eqref{e:ColdDataJump}, but
in most cases it has a negative energy production.
In contrast, Newton's equations always conserve mass, momentum
\emph{and} energy, so the continuum limit of FPU chains beyond a shock
cannot be described in terms of the p-system. Instead, the FPU chain
produces a dispersive shock with strong microscopic
oscillations that take the form of modulated wave trains, 
compare Observations \ref{o:dispShock1}
and \ref{o:dispShock2} below.
\par
Heuristically, the onset of dispersive shocks is a consequence of
energy conservation and can be interpreted as \emph{Hamiltonian
  self-thermalisation}: When the shock is formed the p-system predicts
some macroscopic excess energy which no longer can be stored in cold
motion. On the atomistic scale this excess energy is transferred into
modulated wave trains and appears as \emph{internal} or \emph{thermal}
energy on the macroscopic scale. More precisely, 
although wave trains do not provide a thermalization in the usual `chaotic' 
sense, their macroscopic dynamics is governed by thermodynamically consistent 
field equations, see \cite{DHM06,AMSMSP:DHR}. 
%
\subparagraph{Riemann solvers for the p-system}%
%
The p-system is hyperbolic if $\Phi$ is convex and genuinely
nonlinear for $\Phi^{\prime\prime\prime}\neq0$, thus turning points
of $\Phi^\prime$ correspond to states in which the system is
linearly degenerate. The Riemann problem for strictly convex or
concave flux $\Phi'$ can therefore be described by the
\emph{classical solver}, which is based on the classical Lax theory
from \cite{Lax57} and involves only rarefaction wave and Lax shocks.
\par
This classical Riemann solver is built from the following curves,
where $-$ corresponds to left moving $1$-waves and $+$ to right
moving $2$-waves. The \emph{rarefaction wave sets}
$\calR_\pm[u_{\lst}]$ contain all right states
$u_{\rst}=\pair{r_{\rst}}{v_{\rst}}$ that can be reached from a
given left state $u_{\lst}=\pair{r_{\lst}}{v_{\lst}}$ with a single
$1$- or $2$-rarefaction wave. The \emph{shock wave sets}
$\calS_\pm[u_{\lst}]$ consist of all possible right states $u_\rst$
that can be reached by a single Lax $1$- or $2$-shock. The sets
$\calW_\pm[u_{\lst}]=\calR_\pm[u_{\lst}]\cup\calS_\pm[u_{\lst}]$
form $C^2$-smooth curves through $u_{\lst}$, and we denote
$\calW[u_{\lst}] = \calW_-[u_{\lst}] \cup \calW_+[u_{\lst}]$. The
solution to the Riemann problem with given left and right states
$u_\lst$ and $u_\rst$ consists of the two elementary waves that
connect $u_\lst$ to $u_\mst$, and $u_\mst$ to $u_\rst$ (one of these
may be trivial), where the intermediate state is uniquely determined
by $u_\mst\in\calW_-[u_\lst]$ and $u_\rst\in\calW_+[u_\mst]$.
\par
In Appendix~\ref{app:p-riemann} we give more details about the
classical-solver for the p-system. If $\Phi^\prime$ has turning
points the wave sets of the classical solver must be modified,
and this gives rise to non-classical solvers which involve various
types of composite waves, see \S\ref{s:nonclass}.
%
\subparagraph{Conservative shocks in the p-system}
%
%
In contrast to Lax shocks, \emph{conservative shocks} in the
p-system balance mass, momentum \eqref{e:ColdDataJump} \emph{and}
energy \eqref{e:ColdDataEnJump}. This gives rise to the system of
nonlinear equations
\begin{align}
\label{e:consJumpData}
\begin{split}
c_\rmrh \jump{r} +\jump{v} =0,\qquad c_\rmrh \jump{v}
+\jump{\Phi'(r)} =0,\qquad c_\rmrh \jump{\tfrac{1}{2}v^2+\Phi\at{r}}
+\jump{v\Phi'(r)} =0,
\end{split}
\end{align}
for the five parameters $r_\lst$, $v_\lst$, $r_\rst$, $v_\rst$,
$c_\rmrh$.  According to Appendix \ref{app:p-riemann} a conservative
shock is called \emph{supersonic} if
\begin{math} |c_\rmrh| > \max\{|\lambda_\pm(r_\lst)|,|\lambda_\pm(r_\rst)|\}
\end{math} %
and \emph{subsonic} if
\begin{math} %
|c_\rmrh| < \min\{|\lambda_\pm(r_\lst)|,|\lambda_\pm(r_\rst)|\}
\end{math}. %
It can be easily shown that each conservative shock satisfies
\begin{equation}\label{e:ConsJump}
\calJ(r_\lst,r_\rst):= \jump{\Phi(r)} -
\jump{r}\langle{\Phi'(r)}\rangle =0,\qquad
\langle{\Phi'(r)}\rangle:=\frac{1}{2}(\Phi'(r_\lst)+\Phi'(r_\rst)).
\end{equation}
Conversely, for each solution to \eqref{e:ConsJump} there exist both
a corresponding conservative $1$-shock and $2$-shock. These shocks
are unique up to Galilean transformations, and differ only in
$\sgn\jump{v}=\sgn\,{c_\rmrh}$.
We analyse the set of conservative shocks in the p-system in more
detail in \S\ref{s:mathconsshock}.

%
%
%
\subparagraph{Macroscopic limit and Young measures}
%
About the hyperbolic continuum limit of FPU chains in the presence of
strong microscopic oscillations little is known rigorously. This is
the reason why, for a large part, we have to rely on numerical
observations. The main difficulty lies in the control of oscillations
that lead to measure-valued solutions on the macroscopic
scale. Heuristically, such measure-valued solutions are governed by
extended p-systems, but a rigorous derivation of such extensions could
be treated in a satisfactory manner only for integrable systems so
far; notably the harmonic chain \cite{DHM06,Mie06,Mac02,Mac04}, the
hard sphere model \cite{Her05}, and the Toda chain \cite{DKKZ96,DM98}.
\par
Nevertheless, some insight into the macroscopic evolution of
microscopic oscillations can be gained from the theory of Young
measures. Here it is supposed that the atomic data generate a family
of probability distributions $\mu\triple{\MaTime}{\MaLagr}{\dint{Q}}$,
where $\pair{\MaTime}{\MaLagr}$ is a point in the macroscopic
space-time and $Q=\pair{r}{v}$ denotes a point in the microscopic
phase space of distances and velocities. Note that the atomic data are
oscillation free in the vicinity of a point $\pair{\MaTime}{\MaLagr}$
if and only if the measure $\mu\triple{\MaTime}{\MaLagr}{\dint{Q}}$ is
a delta distribution with respect to the $Q$ variable.
\par
On the one hand, Young measures provide an elegant framework to
investigate oscillatory numerical data that we used to interpret our
simulations. For given $\pair{\MaTime}{\MaLagr}$ the measure
$\mu\triple{\MaTime}{\MaLagr}{\dint{Q}}$ can be approximated by means
of \emph{mesoscopic} space-time windows and provides local mean values
of atomic observables as well as statistical information about the
microscopic oscillations, see \cite{DHM06,DH08}.
\par
On the other hand, Young measures are useful for analytical
considerations because the solutions to \eqref{e:NewEqn1} with
$N\to\infty$ are compact in the sense of Young-measures provided that
the initial data are of order $1$. Extracting
convergent subsequences, one can then prove as in \cite{Her05} that
\emph{every} limit measure must be a \emph{weak solution} to the
following \emph{macroscopic} conservation laws of mass, momentum and
energy
\begin{align}%
\label{e:YMMacroConsLaws} %
\notag%
\MaTimeDer{}\meas{r}-\MaLagrDer{}\meas{v}&=0,\\
\MaTimeDer{}\meas{v}-\MaLagrDer{}\meas{\Phi^\prime\at{r}}&=0,\\
\notag%
\MaTimeDer{}\meas{\textstyle{\frac{1}{2}}v^2+\Phi\at{r}}-%
\MaLagrDer{}\meas{v\,\Phi^\prime\at{r}}&=0.
\end{align}
Here $\meas{\Psi}\pair{\MaTime}{\MaLagr}$ is the \emph{local mean
value} of the observable $\Psi=\Psi\pair{r}{v}$, that is
\begin{align*}
\meas{\Psi}\pair{\MaTime}{\MaLagr}&=%
\int_{\R^2}\psi\at{Q}\,\mu\triple{\MaTime}{\MaLagr}{\rmd Q}.
\end{align*}
System \eqref{e:YMMacroConsLaws} provides non-trivial
information about the macroscopic dynamics of FPU chains. For
arbitrary oscillations, however, we can not express the fluxes in
terms of the densities, and hence \eqref{e:YMMacroConsLaws} does not
determine the macroscopic evolution completely.
\bigpar
As an important consequence of \eqref{e:YMMacroConsLaws} we can
characterise conservative shocks  in FPU chains. To this end suppose
that for all points $\pair{\MaTime}{\MaLagr}$ in a sufficiently
small region the measure $\mu\triple{\MaTime}{\MaLagr}{\dint{Q}}$
depends only on $c=\MaLagr/\MaTime$ and is a delta-distribution
$\delta_{u_\lst}\at{\dint{Q}}$ for $c\leq{c}_\rmrh$ and
$\delta_{u_\rst}\at{\dint{Q}}$ for $c\geq{c}_\rmrh$, for some
$c_\rmrh$. This is exactly what we observe in the numerical FPU
simulations for a conservative shock connecting $u_\lst$ to
$u_\rst$, see Figure~\ref{f:intro_rp_supersonic}. In this case
\eqref{e:YMMacroConsLaws} reduces to the three independent jump
conditions that determine a conservative shock in the p-system,
compare \eqref{e:consJumpData}. In other words, FPU chains allow for
waves that are close to a jump discontinuity only if there exists a
corresponding conservative shock in the p-system.
%
%
\subparagraph{Travelling waves and modulation theory}\label{s:tw}
%
%
Careful investigations of numerical experiments as described in
\cite{DH08} reveal that also for non-integrable cases the oscillations
in a dispersive shock of FPU chains can be described by modulated
travelling waves. That means for each $\pair{\MaTime}{\MaLagr}$ in the
oscillatory region the measure
$\mu\triple{\MaTime}{\MaLagr}{\dint{Q}}$ is generated by a travelling
wave, whose parameters are slowly varying as they depend only on the
macroscopic coordinates $\MaTime$ and $\MaLagr$.  This observation is
in accordance with the fact that the support of
$\mu\triple{\MaTime}{\MaLagr}{\dint{Q}}$ is contained in a closed
curve, compare the density plots in Figures~\ref{f:consShock.B} and
\ref{f:fpu_solver}, which show the superposition of several of these
curves.
\bigpar
Travelling waves with constant speed $c$ are exact solutions to the
infinite chain (\ref{e:NewEqn1}) that depend on a single phase
variable $\phi=k\alpha+\omega{t}$ via $x_\alpha(t)=x(\phi)$. Here
$k$ and $\omega$ are generalized wave number and frequency,
respectively, and $c=-\omega/k$ is the phase velocity.  In terms of
atomic distances and velocities travelling waves can be written as
$r_\MiLagr\at\MiTime=R\at{\phi}$ and
$v_\MiLagr\at\MiTime=V\at{\phi}$, where the profile functions $R$
and $V$ solve the advance-delay differential equations
\begin{align*}%
\notag
c\partial_{\phi}R(\phi)=V(\phi+1)-V(\phi)\,\qquad
c\partial_{\phi}V(\phi)= \Phi^\prime\at{R\at{\phi}}-
\Phi^\prime\at{R\at{\phi-1}}.
\end{align*}
In our context relevant travelling waves are \emph{wave trains}, for
which both $R$ and $V$ are periodic functions, \emph{solitons}
(solitary waves), which limit to the same background state as
$\phi\to\pm\infty$, and \emph{fronts}, which connect to different
constant background states for $\phi\to-\infty$ and $\phi\to\infty$.
\par
Wave trains exist for all convex potentials $\Phi$, see
\cite{FV99,DHM06,Her09}. They depend on four parameters and provide
the building blocks for modulation theory, which describes the
macroscopic evolution of a modulated travelling wave. This evolution
is governed by a system of four nonlinear conservation laws, which
one usually refers to as Whitham's modulation equations, and a
dispersive shock is just a rarefaction wave of this system. The
Whitham equations can be regarded as an extension of
\eqref{e:YMMacroConsLaws}, where modulation theory also provides a
complete set of constitutive relations which depend on $\Phi$ via
the four parameter-family of wave trains. We refer to \cite{Whi74}
for the general background, and \cite{FV99,DHM06} for the
application to FPU chains.
\par%
The existence of solitons for super-quadratic potentials is proven in
broad generality in \cite{FW94,SW97}, and \cite{PP00,Her09} show that
wave trains limit to solitons as the wave number tends to zero, see
also \cite{Pan05}. Solitons are important in our context as they
appear at the shock front of a dispersive shock where the amplitude of
the oscillations is maximal. Generically solitons travel faster than
the sound speed, and converge exponentially to the background state
along distinct directions as curves in the distance-velocity plane,
compare Figure~\ref{f:shock}(c) and \cite{Ioo00}. Near a turning point
of the flux, solitons can travel with the sound speed $\sqrt{\Phi''}$
evaluated at the background state. Such solitons converge
algebraically to the background state as $\phi\to\pm\infty$ and along
the same line in the distance-velocity plane forming a cusp, see
Figure~\ref{f:consShock.B}(c).
\par%
Concerning fronts it has been proven in~\cite{Ioo00} that fronts
bifurcate from turning points $r_*$ of $\Phi^\prime$ with
$\Phi^{(4)}(r_*)\neq 0$ if and only if $\Phi^{(4)}(r_*)<0$, i.e.,
convex-concave flux. These fronts travel faster than the sound
speeds of left and right states, have monotone profiles and converge
to the endstates exponentially, compare also \cite{HR09}. In our
context fronts appear as conservative shocks.
%
%

%
%
\section{Numerical simulations of FPU-Riemann
problems}\label{s:num}
%

%
All simulations in this paper describe Riemann problems with cold
initial data. That means for given \emph{left state}
$u_\lst=\pair{r_\lst}{v_\lst}$ and \emph{right state}
$u_\rst=\pair{r_\rst}{v_\rst}$ we initialise the atomic distances
and velocities by
\begin{align*}
\pair{r_\MiLagr\at{0}}{v_\MiLagr\at{0}}&=\left\{%
\begin{array}{lcl}
\pair{r_\lst}{v_\lst}&&\text{for }\eps\MiLagr\leq \bal_*,\\
\pair{r_\rst}{v_\rst}&&\text{for }\eps\MiLagr > \bal_*,
\end{array}%
\right.%
\end{align*}
where $\eps=1/N$ is the scaling parameter and $\bal_*$ denotes
the macroscopic position of the initial jump. We impose the
boundary conditions
\begin{align*}
v_{N+1}\at\MiTime=v_{N}\at\MiTime,\quad{r_0}\at\MiTime={r_1}\at\MiTime,
\end{align*}
so \eqref{e:NewEqn2} becomes a closed system for the $2N$ unknowns
$r_1\tdots{r_N}$ and $v_1\tdots{v_N}$. These conditions are
appropriate since we start with piecewise constant initial data, and
stop the simulation before any macroscopic wave has reached the
boundary. For the numerical integration of \eqref{e:NewEqn1} we use
the Verlet-scheme, which is a symplectic and explicit integrator of
second order \cite{SSY97,HLW02}. The microscopic time step size
$\triangle\MiTime$ is independent of $N$ and small compared with the
smallest inverse frequency of the linearised chain.
%
%
%
\subsection{Self-similar structure of solutions}
%
%
Typical examples for the numerical outcome of an atomistic Riemann
problem are given in Figures~\ref{Fig:Intro.Data1} and
\ref{Fig:Intro.Data2}, where we plot snapshots of the atomic
distances and velocities against the scaled particle index
$\MaLagr=\MiLagr/N$. The used potential
\begin{align}
\label{Fig:Intro.Data1.Pot}%
\Phi(r)=\exp\at{1-r}-(1-r)+\tfrac{1}{40}\at{r-1}^4
\end{align}
is a modified Toda-potential with strictly convex flux $\Phi^\prime$
and the initial data are given by
\begin{align}
\label{Fig:Intro.Data1.IV}%
\bal_*=0.6,\quad\pair{r_\lst}{v_\lst}=\pair{0}{0},\quad
\pair{r_\rst}{v_\rst}=\pair{0}{1}&.
\end{align}
In Figure~\ref{Fig:Intro.Data1} we fix $N$ and plot the solution for
$\MaTime=0$, and $\MaTime=0.15$, and $\MaTime=0.3$, whereas
Figure~\ref{Fig:Intro.Data2} shows the numerical results for
increasing $N$ at the same macroscopic time $\MaTime=0.3$. Recall that
according to the hyperbolic scaling the microscopic time is always
proportional to $N$.
\par
The simulation indicate that the atomistic solutions indeed converge
on the macroscopic scale to some limiting Young measure which is
self-similar in $\MaTime$ and $\MaLagr-\MaLagr_*$. In the limit
$N\to\infty$ we predict that the solution consist of a cold
rarefaction wave and a dispersive shock which are separated by a
cold intermediate state. In the cold regions the atomic data can be
expected to converge to a macroscopic function, so in each point
$\pair{\MaTime}{\MaLagr}$ the Young measure is a delta distribution.
In the dispersive shock, however, this measure is nontrivial but the
envelopes (and likewise the local mean values) still converge to
functions.
\bigpar
\begin{figure}[t!]%
  \centering{%
  \includegraphics[width=\mhpicDwidth, draft=\figdraft]%
  {\figfile{intro_rp/intro_dist_ini}}%
  \includegraphics[width=\mhpicDwidth, draft=\figdraft]%
  {\figfile{intro_rp/intro_dist_mi_21}}%
  \includegraphics[width=\mhpicDwidth, draft=\figdraft]%
  {\figfile{intro_rp/intro_dist_mi_22}}%
  \\\vskip\mhpicDvsep%
  \includegraphics[width=\mhpicDwidth, draft=\figdraft]%
  {\figfile{intro_rp/intro_vel_ini}}%
  \includegraphics[width=\mhpicDwidth, draft=\figdraft]%
  {\figfile{intro_rp/intro_vel_mi_21}}%
  \includegraphics[width=\mhpicDwidth, draft=\figdraft]%
  {\figfile{intro_rp/intro_vel_mi_22}}%
  }%
  \caption{
  Numerical results for potential
  \eqref{Fig:Intro.Data1.Pot} and Riemann initial data
  \eqref{Fig:Intro.Data1.IV}: Snapshots of atomic distances and
  velocity against the macroscopic particle index $\MaLagr$
  for several macroscopic times and $N=4000$.  On the macroscopic
  scale the solutions becomes self-similar: a left
  moving rarefaction wave and a right moving dispersive shock are
  separated by a cold intermediate state.}
  \label{Fig:Intro.Data1}%
  \centering{%
  $\;$\\
  \includegraphics[width=\mhpicDwidth, draft=\figdraft]%
  {\figfile{intro_rp/intro_dist_mi_12.eps}}%
  \includegraphics[width=\mhpicDwidth, draft=\figdraft]%
  {\figfile{intro_rp/intro_dist_mi_32.eps}}%
  \includegraphics[width=\mhpicDwidth, draft=\figdraft]%
  {\figfile{intro_rp/intro_dist_mi_42.eps}}%
  }%
  \caption{Atomic distances from Figure~\ref{Fig:Intro.Data1}
  for different values of $N$.}
  \label{Fig:Intro.Data2}
\end{figure}%
\bigpar
Figures~\ref{Fig:Intro.Data1} and \ref{Fig:Intro.Data2} provide of
course a merely qualitative confirmation of our interpretation of
the numerical data. A refined quantitative measurement with
different $N$ would be possible but requires much more numerical
effort for the following two reasons.
\par
$(1)$ Since information propagates with infinite speed in the
lattice system \eqref{e:NewEqn1} cold states manifest only in the
limit $N\to\infty$. For finite $N$ we find small
fluctuations everywhere due to the discreteness of $\alpha$.
Heuristically, we expect the amplitude of the fluctuations to decay
exponentially with $N$ outside the space-time cone spanned by the
fastest macroscopic speeds, but only algebraically with $N$ inside
this cone. This expectation is supported by rigorous results for the
Toda chain and the harmonic chain, see \cite{Kam93,MP09}, and
implies that the cold intermediate state has superimposed
fluctuations with amplitude $1/\sqrt{N}$. An accurate measurement of
intermediate states and wave speeds inside the above cone therefore
requires simulations with very large $N$.
\par
$(2)$ Due to the oscillatory nature it is notoriously difficult to
compare quantitatively the dispersive shocks for different values of
$N$ or $\MaTime$: Accurate values for the evolution of the envelopes
are very hard to measure numerically, and the precise values of
averaged quantities such as local mean values or numerical
distribution functions depend for finite $N$ on the details of the
implemented averaging algorithm.
\par
In this paper we focus on the qualitative properties of FPU Riemann
solutions. We aim to understand the macroscopic selection and
composition rules for elementary waves and how turning points of
$\Phi$ effect the qualitative structure of Riemann solutions. In
particular, we do not intend to measure numerical convergence rates
or to predict the wave parameters quantitatively.
%
%
\subsection{Elementary waves}
\begin{figure}
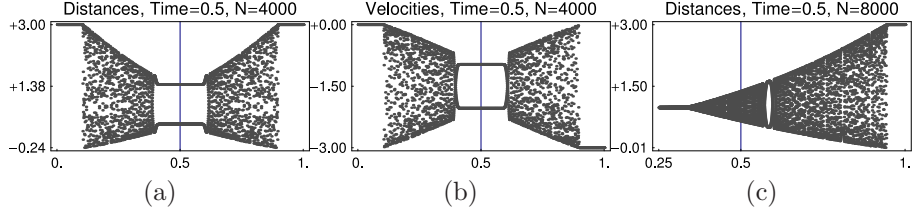

\centering{%
\setlength{\tabcolsep}{0cm}%
\begin{tabular}{ccc}%
\includegraphics[width=\mhpicDwidth, draft=\figdraft]%
{\figfile{toda_piston_and_shock/tdp_dist_mi.eps}}%
&%
\includegraphics[width=\mhpicDwidth, draft=\figdraft]%
{\figfile{toda_piston_and_shock/tdp_vel_mi.eps}}%
&%
\includegraphics[width=\mhpicDwidth, draft=\figdraft]%
{\figfile{toda_piston_and_shock/tds_dist_mi.eps}}%
\\(a) & (b) & (c)%
\end{tabular}%
}%
\caption{%
      Snapshots for the Toda chain
      with initial jump at $\MaLagr_\ast=0.5$ (see vertical lines).
      (a), (b) Distances and velocities for the classical `supercritical shock
      problem' generating a permanent thermalization via binary
      oscillations;
      (c) Single dispersive shock where front and back
      counter-propagate; the asymptotic state at $\MaLagr_\ast$ is a wave train with wave
      number $\sim 0.47$.
      }%
\label{f:thermalize}%
\end{figure}

%
%
%
The following key observation about solutions to FPU-Riemann
problems reflects the hyperbolic and modulation nature of the limit
$\varepsilon\to 0$, but has not yet been proven rigorously.
\begin{observation}
  For strictly monotone and nonlinear flux $\Phi^\prime$ where
  $\Phi'''$ has at most one root in the range of the solution we
  observe the following.
\enum
\item
The macroscopic dynamics for cold Riemann data is self-similar and
hence reducible to the macroscopic velocity variable
$c=\at{\bal-\bal_*}/\bt$.
\item
The arising measure at each $(\bal,\bt)$ is either a point measure
or supported on a closed curve that is generated by the distances
and velocities of a wave train profile. Therefore, we can describe
the macroscopic limit by a family of modulated wave trains
parameterized by $c$.
\item
The macroscopic solution to each cold Riemann problem consist of a
finite number of self-similar waves. These elementary waves are
\begin{enumerate}
\item cold rarefaction waves,
\item dispersive shock fans connecting two cold states,
\item energy conserving jumps between two cold states (only for flux
  with turning point),
\item dispersive shock fans connecting a cold state and a constant
wave train; these shocks always come as a counter-propagating pair.
\end{enumerate}
\end{list}
\end{observation}
%
%
%
%
\subparagraph{Structure of dispersive shocks}
%
Our basic numerical observations concerning dispersive shocks are as
follows, see also Figure~\ref{f:shock}. As mentioned in
\S\ref{chIntro}, dispersive shocks have been studied for certain
potentials and in other contexts, but we have not found the following
explicitly mentioned for FPU.
\begin{observation}\label{o:dispShock1}
In the numerical simulation of cold Riemann problems
\emph{dispersive shocks} appear with oscillatory atomic data between
two constant states, see Figure \ref{f:shock}. Within a dispersive
shock the self-similarity variable
$c=\at{\MaLagr-\MaLagr_*}/\MaTime$ ranges between the \emph{shock
back velocity} $c_\rmb$ and the \emph{shock front velocity}
$c_\rmf$. The atomic oscillations have monotone envelopes with
maximal amplitude at the front and vanishing at the back. There
exist dispersive $1$-shocks with $c_\rmf<0$ and $c_\rmf<c_\rmb$ as
well as dispersive $2$-shocks with $c_\rmf>0$ and $c_\rmf>c_\rmb$.
\end{observation}
Moreover, our numerical results suggest the following fine-structure
of the oscillations within a dispersive shock.
\begin{observation}\label{o:dispShock2}
  A dispersive shock consists of a one-parameter family of wave trains
  parameterised by $c=\at{\MaLagr-\MaLagr_*}/\MaTime$. Within the
  dispersive shock the parameter modulation is smooth (and hence
  follows a rarefaction wave of Whitham's modulation equations).
  Dispersive $2$-shocks have the following properties (1-shocks
  accordingly due to symmetry).  \enum
\item The measure $\mu\at{c_\rmb}$ at the back of the shock reduces to
  the point measure generated by the constant left state $u_{\lst}$,
  and the local mean values of atomic distances and velocities
  smoothly connect to $u_{\lst}$.
\item The measure $\mu\at{c_\rmf}$ at the front of the shock converges
  to a \emph{soliton} with background state $u_{\rst}$, and the
  local mean values are continuous but not differentiable at $c_\rmf$.
\item The family of curves $\mathrm{supp}(\mu(c))$ with
  $c_\rmb<c<c_\rmf$ is nested.
\item The dispersive shock is \emph{compressive} in the sense that
  both positive characteristic speeds $\lambda_\pm$ of the p-system
  point into the fan, i.e. $\lambda_+\at{u_{\lst}}>c_\rmb$ and
  $\lambda_+\at{u_{\rst}}<c_\rmf$.
\item The Rankine-Hugeniot velocity of the jump lies strictly between
  $c_\rmb$ and $c_\rmf$.
\end{list}
\end{observation}
For sufficiently small jump heights
$\abs{r_\rst-r_\lst}+\abs{v_\rst-v_\lst}$ the shock back and front
move in the same direction, this means we have either
$c_\rmf<c_\rmb<0$ or $c_\rmf>c_\rmb>0$. The classical `piston
problem', however, shows that the situation is more subtle for large
jump heights. The initial data in this problem describe an
evenly spaced chain with positive left velocities and negative right
velocities. It has been observed in \cite{HS78} that sufficiently
large (`supercritical') jumps in the velocity generate a transition
from a pair of counter-propagating dispersive shocks with cold
intermediate state to incomplete dispersive shocks whose `backs' are
constant binary oscillations which replace the cold intermediate
state. This phenomenon is related to dispersive shock fans with
counter-propagating back and front, see
Figure~\ref{f:thermalize}(c). More precisely, passing the critical
jump height from below the shock back velocities change their sign
and the oscillatory intermediate state for supercritical data
results from the interaction of two dispersive shocks. Analysing the
back velocity of a dispersive shock might be a fruitful approach to
determine the critical jump height in the piston problem.
%
\subparagraph{Conservative shocks in FPU chains}
%
%
As discussed in \S\ref{chIntro}, numerical simulations indicate
that conservative shocks appear in FPU Riemann problems only if they
are supersonic. For illustration we consider the two quintic
potentials
\begin{eqnarray}
\Phi(r+2) &=& r^2 - \frac{r^3}{6} - \frac{r^4}{24} + \frac{r^5}{120}%
\label{e:fastpot}\\
\Phi(r+2) &=& r^2 - \frac{r^3}{6} + \frac{r^4}{24} +
\frac{r^5}{120}%
\label{e:slowpot}.
\end{eqnarray}
We plot the set of conservative shocks for these potentials in
Figure~\ref{f:consdata}. Note that due to
Theorem~\ref{t:consData}\eqref{i:conn} each of these sets consists
of the diagonal and a closed curve crossing the diagonal at
the turning points.
\bigpar
The flux $\Phi'$ for \eqref{e:fastpot} has a turning point
$r_*\approx 1.3$ with $\Phi^{(4)}(r_*)<0$ (convex-concave), while
the flux for \eqref{e:slowpot} has a turning point $r_*\approx 2.7$
with $\Phi^{(4)}(r_*)>0$ (concave-convex). Both potentials are
convex in a neighborhood of $r^*$, and the other turning points are
outside the range of simulation.

\begin{figure}[h!]
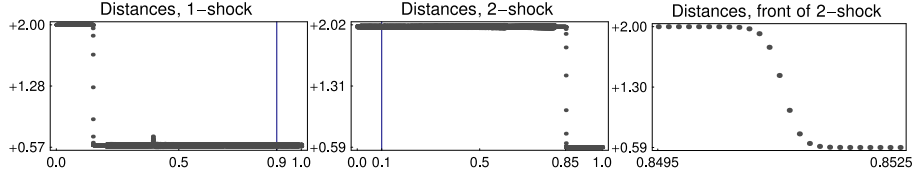

\setlength{\tabcolsep}{0cm}%
\centering{
\includegraphics[width=\mhpicDwidth, draft=\figdraft]%
{\figfile{qp1_conservative_shocks/qp1_shock_1_dist.eps}}%
\includegraphics[width=\mhpicDwidth, draft=\figdraft]%
{\figfile{qp1_conservative_shocks/qp1_shock_2_dist.eps}}%
\includegraphics[width=\mhpicDwidth, draft=\figdraft]%
{\figfile{qp1_conservative_shocks/qp1_shock_2_dist_zoom.eps}}%
}%
\caption{%
  Supersonic conservative $1$- and $2$-shock for potential
  \eqref{e:fastpot} with  $N=8000$, $\MaTime=0.5$ and
  $\bal_\ast$ indicated by the vertical lines.  }
\label{f:consShock.A} %
\end{figure}
\par%
  Potential \eqref{e:fastpot} allows for instance for the two
  supersonic conservative shocks
\begin{align*}
r_\lst=2,\quad r_\rst\approx0.59,\quad c_\rmrh=\pm 1.50,\quad
v_\rst-v_\lst\approx\pm2.11,\quad \lambda_\pm(r_\lst)\approx \pm
1.41,\quad \lambda_\pm(r_\rst)\approx \pm0.89.
\end{align*}
In Figure \ref{f:consShock.A} we plot the solution to the
corresponding FPU-Riemann problem (with $v_\lst=0$) and conclude
that both the supersonic $1$-shock and the supersonic $2$-shock are
captured by the atomic chain.
\par
The conservative shocks for potential \eqref{e:slowpot} that range
over $r_*=2.7$, however, are subsonic. The FPU solution to
\begin{align}
r_\lst=4,\quad{}r_\rst\approx1.24,\quad c_\rmrh\approx\pm 1.46,\quad
v_\rst-v_\lst\approx\pm4.03,\quad \lambda_\pm(r_\lst)\approx \pm
1.83,\quad \lambda_\pm(r_\rst)\approx \pm1.73,
\end{align}
is plotted in Figure \ref{f:consShock.B}, and is far from a
conservative shock. Recall that this is accordance with the
non-bifurcation result for subsonic fronts in \cite{Ioo00}. The
solution in Figure \ref{f:consShock.B} consists of a dispersive shock
with \emph{attached} rarefaction wave, that means both waves are not
separated by a constant state. Consequently, the soliton at the front
of the dispersive shock is no longer supersonic but sonic, i.e., its
speed equals the sound speed of the background state. This is
confirmed by the numerical data in Figure \ref{f:consShock.B}(c): the
distribution function near the soliton has the predicted cusp
shape. Compare with the non-degenerate exponentially decaying soliton
in Figure~\ref{f:shock}(c).
%
%
%
\section{Riemann solvers}\label{s:riem}
%
%
%
%
\subsection{Towards an FPU-Riemann solver for the classical case}\label{s:RieSolDispShock}
\begin{figure}
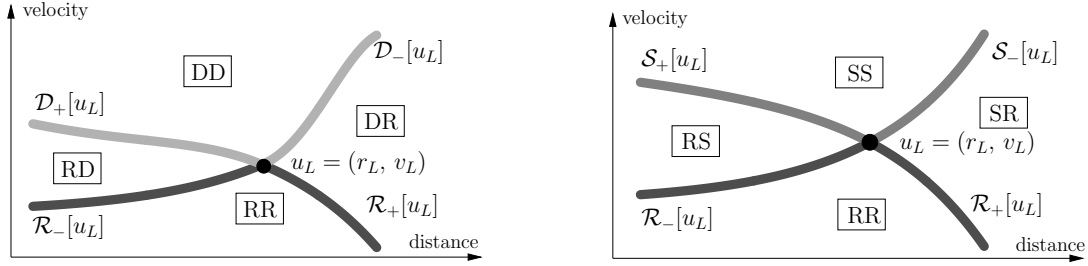

\centering{%
\includegraphics[width=0.4\textwidth, draft=\figdraft]%
{\figfile{xfig/riemann_1/rs_atomistic_convex_flux.eps}}%
\hspace{0.1\textwidth}%
\includegraphics[width=0.4\textwidth, draft=\figdraft]%
{\figfile{xfig/riemann_1/rs_standard_convex_flux.eps}}%
}%
\caption{%
\emph{Left}: Sketch of the FPU Riemann-solver for strictly
convex $\Phi'$. The wave sets $\calW^\micro[u_\lst]$ consist of
rarefaction curves and dispersive shock curves, and decompose the
plane into 4 regions DD, RD, RR, and RD. \emph{Right}: The
corresponding classical solver for the p-system with Lax shocks
instead of dispersive shocks.
}%
\label{fig:atomistic_rs}%
\end{figure}

%
In this section we describe an adaption of the classical p-system
solver that accounts for dispersive shocks in the macroscopic
solutions to FPU-Riemann problems. Based on Observations
\ref{o:dispShock1} and \ref{o:dispShock2} we arrive at the
following conjecture, which is illustrated
Figure~\ref{fig:atomistic_rs}.
\begin{conjecture}
From each state $u_{\lst}$ there emanate two \emph{dispersive shock
curves} $\calD_-[u_{\lst}]$ and $\calD_+[u_{\lst}]$ with the
following properties.
\enum
\item
Each state $u_{\rst}\in\calD_\pm[u_{\lst}]$ can be connected with
$u_{\lst}$ by a single dispersive shock with $\sgn(c_\rmf)=\pm 1$.
\item
  The curves fit smoothly to the corresponding rarefaction curves
  $\calR_\pm[u_{\lst}]$ (so that for small jump heights the dispersive
  and Lax-shock curves almost coincide).
\end{list}
The macroscopic solution to FPU-Riemann problems with cold data and
sufficiently small jump heights can be described by the wave sets
\begin{align*}
\calW^\micro_\pm[u_{\lst}]=\calR_\pm[u_{\lst}]\cup\calD_\pm[u_{\lst}].
\end{align*}
In particular, a FPU-Riemann solution consists of a unique $1$-wave
from $\calW^\micro_-[u_{\lst}]$, an intermediate state $u_{\mst}$,
and a unique $2$-wave $\calW^\micro_+[u_{\mst}]$.  The intermediate
state is cold if either a rarefaction wave occurs or if adjacent
shock backs move away from each other; otherwise it is a wave train.
\end{conjecture}
\begin{figure}[ht!]%
\centering{%
\includegraphics[width=\mhpicDwidth, draft=\figdraft]%
{\figfile{dispersive_shock_curve/dsc_dist_mi_prm_06.eps}}%
\includegraphics[width=\mhpicDwidth, draft=\figdraft]%
{\figfile{dispersive_shock_curve/dsc_dist_mi_prm_14.eps}}%
\includegraphics[width=\mhpicDwidth, draft=\figdraft]%
{\figfile{dispersive_shock_curve/dsc_dist_mi_prm_20.eps}}%
}%
\caption{%
  Three points from the dispersive shock curve $\calD_-[\pair{0}{0}]$
  for potential \eqref{Fig:Intro.Data1.Pot}. }%
\label{fig:dispersive_shock_curve}%
\end{figure}
To illustrate the first part of this conjecture we present
simulations for the modified Toda potential
\eqref{Fig:Intro.Data1.Pot}. For given left state
$u_{\lst}=\pair{0}{0}$ and different values of $r_\rst<0$ we choose
$v_\rst$ such that
$u_{\rst}=\pair{r_\rst}{v_\rst}\in{}\calS_-[u_{\lst}]$, and study
the macroscopic behaviour for the corresponding solutions to
Newton's equations. The results for $\bal_*=0.5$, $\MaTime=0.1$, and
$N=4000$ are depicted in Figure \ref{fig:dispersive_shock_curve}.
For all values of $r_\rst$ we find a dispersive $1$-shock whose
front moves to the left, and an essentially cold intermediate state
$u_{\mst}$ which gives a point in $\calD_-[u_{\lst}]$. In all
simulations there exists a right moving $2$-wave, but this wave has
much smaller amplitudes than the $1$-wave. Hence,
$\calD_-[u_{\lst}]$ and $\calS_-[u_{\lst}]$ are different but close
to each other. Moreover, the simulations indicate the following
behaviour for increasing jump height $r_{\lst}-r_\rst$. The shock
front velocity $c_\rmf$ decreases, whereas $c_\rmb$ increases and
changes its sign at a critical value of $r_\rst$. As mentioned, this
critical value can be viewed as the analogue to the critical
parameter in the piston problem, see Figure~\ref{f:thermalize}.
\par %
\begin{figure}
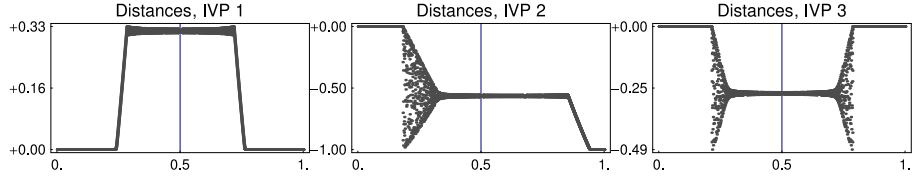

\centering{%
\includegraphics[width=\mhpicDwidth, draft=\figdraft]%
{\figfile{dispersive_riemann_solver/drs_dist_mi_ivp_1.eps}}%
\includegraphics[width=\mhpicDwidth, draft=\figdraft]%
{\figfile{dispersive_riemann_solver/drs_dist_mi_ivp_2.eps}}%
\includegraphics[width=\mhpicDwidth, draft=\figdraft]%
{\figfile{dispersive_riemann_solver/drs_dist_mi_ivp_3.eps}}%
}%
\caption{%
Illustration of the FPU-Riemann solver for potential
\eqref{Fig:Intro.Data1.Pot} and initial data
\eqref{ex:atomistic_solver_inidata}. The examples IVP1--3 correspond
to the regions RR, DR, and DD,  respectively, in Figure
\ref{fig:atomistic_rs}.
}%
\label{fig:dispersive_riemann_solver}
\end{figure}
The resulting FPU-Riemann solver is illustrated in Figure
\ref{fig:dispersive_riemann_solver}, which shows the solutions to
the following initial value problems with $u_\lst=\pair{0}{0}$ and
\begin{align}
\label{ex:atomistic_solver_inidata} %
\begin{array}{lclclclclcl}
1: &u_\rst=\pair{0}{+1},&\quad& 2: &u_\rst=\pair{-1}{0}, &\quad&
3:&u_\rst=\pair{0}{-1}.
\end{array}
\end{align}%
It is important to note that the Lax shock wave sets
$\calS_\pm[u_{\lst}]$ generally differ from the dispersive shock
wave sets $\calD_\pm[u_{\lst}]$. Therefore, replacing
$\calS_\pm[u_{\lst}]$ by $\calD_\pm[u_{\lst}]$ changes the Riemann
solution, i.e., the precise values for the intermediate state and
possibly the waves themselves.
\par
The difference between FPU chain and p-system can be quantified
for the Toda potential \eqref{e:TodaPot}. For the shock piston with
$u_{\lst}=\pair{1}{+2a}$, $u_{\rst}=\pair{1}{-2a}$, and $0<a<1$
(subcritical case) the classical solver for the p-system provides the
intermediate state $u_{\mst}=\pair{r_M}{0}$ with
\begin{align*}
4a^2=\at{r_M-1}\Bat{1-\exp\at{1-r_M}}=\at{r_M-1}^2+\mathcal{O}\bat{\at{r_M-1}^3}
\end{align*}
whereas the results in \cite{Kam93} imply
\begin{align*}
r_M-1=2\ln\at{a+1}\at{a+1}^2=2a+\mathcal{O}\bat{a^2}
\end{align*}
for the FPU-Riemann solver. Both results are different but agree to
leading order in $a$.
%
%
%
\subsection{Riemann solvers for the p-system in the non-classical cases}\label{s:nonclass}
%
%
In this section we consider forces $\Phi'$ with a single turning
point or where the solution ranges over at most one turning point.
For the p-system, the classical Riemann solver cannot be used and
also the above FPU-Riemann solver fails. To prepare the discussion
of FPU-Riemann solvers in these cases, we first describe the
relevant solvers of the p-system.
\par
The building blocks for each non-classical solver are modified wave
sets which replace $\calW[u_\lst]$ from the classical solver.
Specifically, wave curves need to be adapted when intersecting the
line $r=r_\ast$, though these start out near $u_\rst=u_\lst$ as
classical wave sets (rarefaction waves and Lax shocks). The reason
is the change in the sign of $\Phi'''$ at $r_*$, which implies that
the Lax condition \eqref{eqn:LaxCond} or the compatibility condition
$\lambda\at{r_\rst}>\lambda\at{r_\lst}$ is violated. Note that the
wave curves directed away from the turning point $r_\ast$ remain
unchanged.
For a single turning point of the flux the classical wave sets
interact with the turning point as follows.
\begin{remark}
\label{Rem:GeometryWaveSets}
Let $u_\lst$ with $r_\lst\neq{r_\ast}$ be given and suppose
\emph{convex-concave} $\Phi'$\ with $\Phi^{(4)}\at{r_\ast}<0$. Then
the curves $\calR_+[u_\lst]$ and $\calS_-[u_\lst]$ intersect the
line $r=r_\ast$ in the $(r,v)$-plane, whereas $\Phi'''$ does not
change its sign along $\calR_-[u_\lst]$ and $\calS_+[u_\lst]$. The
same holds for \emph{concave-convex} $\Phi'$\ with
$\Phi^{(4)}\at{r_\ast}>0$ if we replace + by $-$.
\end{remark}
\begin{Proof}
We start with convex-concave $\Phi$. For $r_\lst<r_\ast$ we have
$\Phi^{\prime\prime\prime}\at{r_\lst}>0$ and the formulas of
Appendix \ref{app:p-riemann} imply that both $\calR_+[u_\lst]$,
$\calS_-[u_\lst]$ point into direction of increasing $r$, whereas
$r$ decreases along $\calR_-[u_\lst]$ and $\calS_+[u_\lst]$; compare
Figure \ref{fig:atomistic_rs} for an illustration of the wave sets
of $u_\lst$. The same is true for $r_\lst>r_\ast$ as
$\Phi^{\prime\prime\prime}\at{r_\lst}<0$ implies that now $r$
decreases along $\calR_+[u_\lst]$, $\calS_-[u_\lst]$. Finally, the
proof for concave-convex $\Phi^\prime$ is analogous.
\end{Proof}
%
%
\subparagraph{The conservative solver}
%
\begin{figure}[ht!]
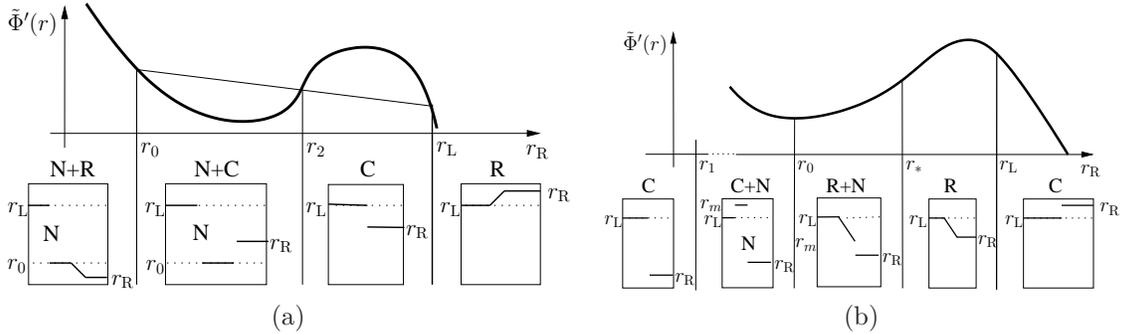
%
\begin{center}%
\begin{tabular}{cc}%
\includegraphics[width=0.47\textwidth, draft=\figdraft]%
{\figfile{xfig/riemann_2/psys_cons_mod_shock}}%
&%
\includegraphics[width=0.43\textwidth, draft=\figdraft]%
{\figfile{xfig/riemann_2/psys_cons_mod_rare}}%
\\(a) & (b)%
\end{tabular}%
\end{center} %
\caption{%
Modified wave sets of the conservative p-system solver: (a)
modified shock curve, (b) modified rarefaction curve;
C$=$compressive Lax shock, N$=$non-classical conservative shock,
R$=$rarefaction wave; $\widetilde\Phi'(r)$ as in the main text. In
the insets we sketch an example for the different composite waves in
each segment of the $r_\rst$-axis.
}
\label{f:ConsSolver.pSys}%
\end{figure}%
The \emph{modified shock curve} $\calS_-^\cons[u_\lst]$ of the
conservative p-system solver is illustrated in Figure
\ref{f:ConsSolver.pSys}(a) for convex-concave $\Phi^\prime$ and
$r_\lst>r_\ast$. This is done by parameterizing
$\calS_-^\cons[u_\lst]$ by $r_\rst$ and plotting $\widetilde\Phi'(r)
= \Phi'(r)+\rho_1 r + \rho_2$ to give an illustrative graph. In
accordance with Remark~\ref{Rem:GeometryWaveSets} the curve
$\calS_-^\cons[u_\lst]$ starts out as $\calS_-[u_\lst]$ for
$r_\rst\lessapprox{} r_\lst$. The wave set modification for large
$r_\lst-r_\rst$ goes via conservative shocks and works as follows.
Along $\calS^\cons_-[u_\lst]$ there exist a unique state
$u_0=\pair{r_0}{v_0}\in\calH_-[u_\lst]$ that can be reached from
$u_\lst$ with a single conservative $1$-shock. This state $u_0$
determines another state $u_2=\pair{r_2}{v_2}$ corresponding to a
Lax shock in $\calS_-[u_\lst]$ such that at $r_2$ the secant slope
from $r_\lst$ to $r_2$ coincides with the slope of the secant from
$r_\lst$ to $r_0$. This reads
\begin{align*}
r_0<r_2<r_\lst\quad\text{and}\quad
\frac{\Phi'\at{r_\lst}-\Phi'\at{r_2}}{r_\lst-r_2}=
\frac{\Phi'\at{r_0}-\Phi'\at{r_2}}{r_0-r_2},
\end{align*}
and implies that the Lax shock connecting $u_\lst$ to $u_2$ and the
conservative shock connecting $u_\lst$ to $u_0$ have the same
Rankine-Hugeniot velocity. Note that this relation and in fact all
conservative shock distance data are the same for $\Phi'$ and
$\widetilde\Phi'$.
\par
\begin{figure}[ht!]
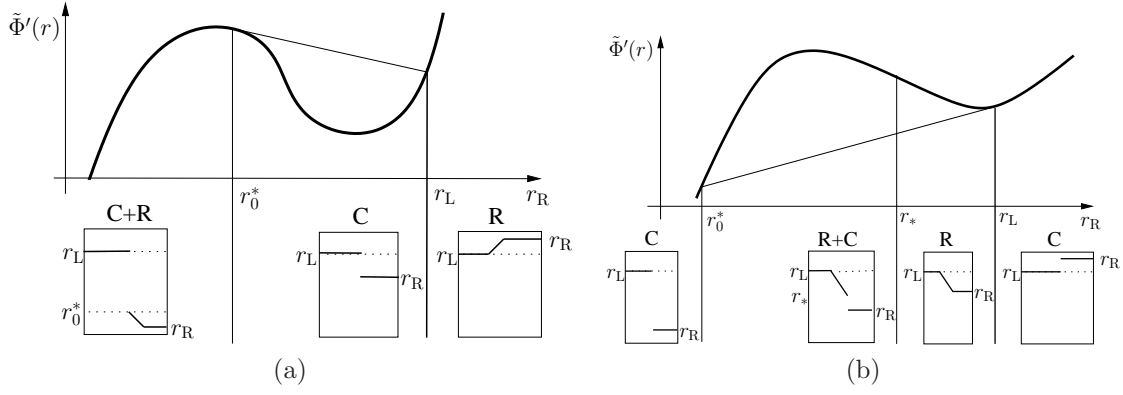
%
\begin{center}%
\begin{tabular}{cc}%
\includegraphics[width=0.47\textwidth, draft=\figdraft]%
{\figfile{xfig/riemann_2/psys_diss_mod_shock}}%
&%
\includegraphics[width=0.43\textwidth, draft=\figdraft]%
{\figfile{xfig/riemann_2/psys_diss_mod_rare}}%
\\(a) & (b)%
\end{tabular}%
\end{center} %
\caption{%
Modified wave sets of the dissipative p-system solver: (a) modified
shock curve, (b) modified rarefaction curve; symbols as in Figure
\ref{f:ConsSolver.pSys}.
}
\label{f:DissSolver.pSys}%
\end{figure}

The \textit{modified rarefaction curve} $\calR_+^\cons[r_\lst]$ is
illustrated in Figure~\ref{f:ConsSolver.pSys}(b) in the same way. The
following lemma, which follows directly from results in \cite{Smo94,
  LeF02}, precisely describes the modified shock and rarefaction
curves (for both convex-concave and concave-convex
$\Phi^\prime)$.
\begin{lemma}\label{l:consRieSol}~
\enum
\item
  If $\Phi'''(r)\neq 0$ for all $r\in I_0:=(r_\rst,r_\lst)$, then the
  solution coincides with the classical solution and consists of at
  most two uniquely chosen rarefaction fans or compressive shocks.%
\item
  Suppose $\Phi'''(r_*)=0$ for a unique $r_*$ and let $u_\lst=(r_\lst,
  v_\lst)$ with $r_\lst>r_*$ be given. Then the following are uniquely
  defined for $r_\rst<r_\lst$. The solution $r_0$ of
  $\calJ(r_\lst,r_0)=0$, the solution $r_m$ of $\calJ(r_\rst,r_m)=0$,
  and the solutions $r_1$ of $|c_\rmrh(r_\lst,r_1)| =
  |c_\rmrh(r_\lst,r_m)|$ and $r_2$ of $|c_\rmrh(r_\lst,r_2)| =
  |c_\rmrh(r_\lst,r_0)|$. It holds that $r_1 < r_0 < r_2 < r_* <
  r_\lst$.
  Note that $r_0, r_2$ depend only on $r_\lst$ whereas $r_m, r_1$
  are functions of $r_\rst$.
\begin{enumerate}
\item[(a)]
  Let $v_0$ be such that $u_0:= (r_0,v_0)\in\calH_\pm[u_\lst]$. A
  right state $u_\rst$ with $r_\rst<r_\lst$ lies in the modified shock
  wave set $\calS^\cons_\pm[u_\lst]$ if $u_\rst\in\calS_\pm[u_\lst]$
  for $r_2<r_\rst$, $u_\rst\in\calS_\pm[u_0]$ for $r_0<r_\rst<r_2$,
  and $u_\rst\in\calR_\pm[u_0]$ for $r_\rst<r_0$.
  %
  %
  By definition of $r_0$, the shock from $u_\lst$ to $u_0$ is
  conservative, and the solution amplitude is discontinuous at $r_2$.
  There always is an intermediate state between rarefaction fan or
  compressive shock and conservative shock.
\item[(b)]
  Let $v_m$ be such that $u_m:= (r_m,v_m)\in\calH_\pm[u_\rst]$.
  A right state $u_\rst$ with $r_\rst<r_\lst$ lies in the modified
  rarefaction wave set $\calR^\cons_\pm[u_\lst]$ if
  $u_\rst\in\calR_\pm[u_\lst]$ for $r_\rst > r_*$,
  $u_m\in\calR_\pm[u_\lst]$ for $r_0<r_\rst<r_*$,
  $u_m\in\calS_\pm[u_\lst]$ for $r_1<r_\rst<r_0$, and
  $u_\rst\in\calH_\mp[u_\lst]$ for $r_\rst<r_1$.
  \par%
  By definition of $r_m$, for $r_1<r_\rst<r_*$ the shock from $u_m$ to
  $u_\rst$ is conservative, and the solution amplitude is
  discontinuous at $r_1$. For $r_0<r_\rst<r_*$ the rarefaction fan is
  attached to the conservative shock, while for $r_1<r_\rst<r_0$ the
  compressive shock is not.
\end{enumerate}
\item If the flux $\Phi'$ has several turning points, then the
  modified wave sets are unchanged as long as $u_\lst$, $u_\rst$ are
  such that only one turning point lies in $[r_1,r_m]$ and $[r_0,r_\lst]$.
\end{list}
\end{lemma}
\begin{Proof}
$(1.)$ Lemma~\ref{l:consShock}(2.) shows that $\calJ\neq 0$ in this
case. Hence, the conservative Riemann solver coincides with the
classical solver \cite{LeF02}. 
For this solver, the p-system is solved uniquely in terms of at most
two rarefaction or shock waves \cite{Smo94}.
$(2.)$ The proof for the shock case is an immediate consequence of
Theorem IV.4.3 (see also Theorem II.4.3) of \cite{LeF02}.
The rarefaction case is not explicitly proven in \cite{LeF02}, but
follows from the exposition, cf. \cite{LeF02} p.163 (see also
Theorem II.5.4). $(3.)$  This follows from the independence of the
solution on $\Phi$ outside this range.
\end{Proof}
Due to this construction, the Riemann solution generically contains
conservative shocks despite the fact that these are of higher
codimension in the space of left and right states: A solution will
consist of three elementary waves instead of two whenever a
conservative shock is possible. Also note that the solution is
non-monotone whenever a compressive and a conservative shock are
connected in a solution.
%
%
\subparagraph{The dissipative solver}
%
%
We refer to the `maximum entropy dissipation' solver in \cite{LeF02}
as the \emph{dissipative solver}. This solver is much simpler than
the conservative solver, and we do not explain it in as much detail.
We plot the modified wave sets which determine the solver in
Figure~\ref{f:DissSolver.pSys} for concave-convex $\Phi^\prime$
and $r_\lst>r_*$. Compared with the conservative solver, the
regions with conservative shocks have shrunk to points. The distance
values where wave sets need to be modified are the turning point
$r_*$ and the value $r_0^*$ where the compressive shock has
\emph{extremal} velocity, i.e., where it coincides with a
characteristic velocity.
%
%
\subsection{Towards an FPU-Riemann solver for the non-classical
supersonic case}
%
%
%
We numerically tested the predictions of the conservative solver
about the modified wave sets by simulations with initial data on the
p-system shock and rarefaction curves. From the classical case we
expect that compressive shocks are replaced by dispersive shocks in
the FPU chain.  Indeed, up to this modification, for convex-concave
potentials the conservative solver qualitatively makes the correct
predictions for FPU chains.
\begin{figure}[ht!]
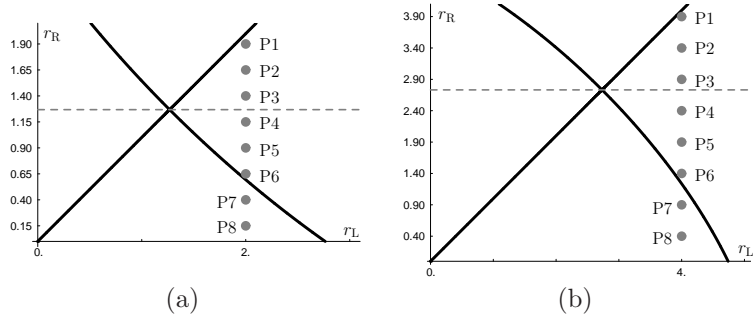

\centering{%
\begin{tabular}{cc}%
\begin{minipage}[c]{0.3\textwidth}%
\includegraphics[width=\textwidth, draft=\figdraft]%
{\figfile{xfig/ivp_data/qp1_ivp_data.eps}}%
\end{minipage}%
&%
\begin{minipage}[c]{0.3\textwidth}%
\includegraphics[width=\textwidth, draft=\figdraft]%
{\figfile{xfig/ivp_data/qp2_ivp_data.eps}}%
\end{minipage}%
\\(a) & (b)%
\end{tabular}%
}%
\caption{%
  Part of conservative shock data for the potentials \eqref{e:fastpot}
  in (a) and \eqref{e:slowpot} in (b). The horizontal lines mark the
  turning point. Bullets
  mark the locations of $(r_\lst,r_\rst)$ for the simulations
  the simulations plotted in Figures~\ref{fig:qp1_shock_curve} and
  \ref{fig:qp1_rarefaction_curve} for (a), as well as
  Figures~\ref{fig:qp2_shock_curve} and
  \ref{fig:qp2_rarefaction_curve} for (b).  } %
\label{f:consdata}%
\end{figure}
\bigpar%
\begin{figure}
\centering{%
\includegraphics[width=\mhpicDwidth, draft=\figdraft]%
{\figfile{qp1_shock_curve/qp1_sc_dist_mi_S2.eps}}%
\includegraphics[width=\mhpicDwidth, draft=\figdraft]%
{\figfile{qp1_shock_curve/qp1_sc_dist_mi_S3.eps}}%
\includegraphics[width=\mhpicDwidth, draft=\figdraft]%
{\figfile{qp1_shock_curve/qp1_sc_dist_mi_S4.eps}}%
\vskip\mhpicDvsep%
\includegraphics[width=\mhpicDwidth, draft=\figdraft]%
{\figfile{qp1_shock_curve/qp1_sc_dist_mi_S5.eps}}%
\includegraphics[width=\mhpicDwidth, draft=\figdraft]%
{\figfile{qp1_shock_curve/qp1_sc_dist_mi_S6.eps}}%
\includegraphics[width=\mhpicDwidth, draft=\figdraft]%
{\figfile{qp1_shock_curve/qp1_sc_dist_mi_S7.eps}}%
}%
\caption{%
  Simulations for data on the supersonic $1$-shock curve for
  potential \eqref{e:fastpot} and $u_\lst=(2,0)$,
  $\MaLagr_\ast=0.9$, $\MaTime=0.5$,
  $N=2000$. The values for $r_\rst$ are those in
  Figure~\ref{f:consdata}(a), so `IVP Sx' crosses the turning point
  for increasing $x$.
}%
\label{fig:qp1_shock_curve}%
\quad\newline
\centering{%
\includegraphics[width=\mhpicDwidth, draft=\figdraft]%
{\figfile{qp1_rarefaction_curve/qp1_rc_dist_mi_R1.eps}}%
\includegraphics[width=\mhpicDwidth, draft=\figdraft]%
{\figfile{qp1_rarefaction_curve/qp1_rc_dist_mi_R3.eps}}%
\includegraphics[width=\mhpicDwidth, draft=\figdraft]%
{\figfile{qp1_rarefaction_curve/qp1_rc_dist_mi_R5.eps}}%
\vskip\mhpicDvsep%
\includegraphics[width=\mhpicDwidth, draft=\figdraft]%
{\figfile{qp1_rarefaction_curve/qp1_rc_dist_mi_R6.eps}}%
\includegraphics[width=\mhpicDwidth, draft=\figdraft]%
{\figfile{qp1_rarefaction_curve/qp1_rc_dist_mi_R7.eps}}%
\includegraphics[width=\mhpicDwidth, draft=\figdraft]%
{\figfile{qp1_rarefaction_curve/qp1_rc_dist_mi_R8.eps}}%
}%
\caption{%
  Simulations for data on the $2$-rarefaction curve for
  potential \eqref{e:fastpot} analogous to Figure
  \ref{fig:qp1_shock_curve}.
}%
\label{fig:qp1_rarefaction_curve}%
\end{figure}
In order to illustrate the structure of the modified wave sets for
\eqref{e:fastpot} we proceed as follows. We fix the left state
$u_\lst=\pair{2}{0}$, and for the points $Px$ marked in Figure
\ref{f:consdata}(a) we solve two Riemann problems denoted by $Sx$ and
$Rx$. The value for ${r_\rst}$ is determined by $Px$, whereas
${v_\rst}$ is chosen such that, for the p-system, $Sx$ and $Rx$
correspond to a single $1$-shock and $2$-rarefaction wave,
respectively. The numerical results for the atomic chain are plotted
in Figures~\ref{fig:qp1_shock_curve}
and~\ref{fig:qp1_rarefaction_curve}.
\par%
Neglecting small waves and fluctuations (caused by the computational
boundary and the positivity of $\eps$) the simulations provide
numerical evidence that the macroscopic limit can indeed be
described by a modified conservative solver. The solutions for $S1$
to $S8$ correspond to those in Figure \ref{f:ConsSolver.pSys}(a)
when replacing compressive by dispersive shocks, and we use the
inset titles C, N+C, N+R in the following. For small jump, i.e.
$r_\rst\lessapprox{}r_\lst$, the chain produces dispersive
$1$-shocks with amplitudes proportional to the jump height ($S2, S3
\cong$ C), and there exists a critical point $r_\rst=\hat{r}_2$ at
which the unique supersonic conservative shock nucleates. For
$r_\rst<\hat{r}_2$ the conservative shock persists whereas the
dispersive shock shrinks ($S4, S5, S6 \cong$ N+C) and transforms
into a rarefaction wave ($S7 \cong$ N+R). Note that, in contrast to
the conservative p-system solver, the nucleation of the conservative
shock is a continuous transition in the envelope since the
dispersive shock extends to the nucleating intermediate state. In
the same way the solutions of $R1$ to $R8$ transform to those of
Figure \ref{f:ConsSolver.pSys}(b), again using inset titles: $R1, R3
\cong$ R, $R5, R6 \cong$ R+N (with $R6$ near N only), $R7, R8 \cong$
C+N.
\bigpar%
We summarize these modifications for the non-classical
supersonic case in the following conjecture.

\begin{conjecture}
For convex-concave flux Lemma~\ref{l:consRieSol} holds under the
following modifications and thereby defines the modified shock wave
sets $\calS_\pm^\fpu[u_\lst]$ and the modified rarefaction wave sets
$\calR_\pm^\fpu[u_\lst]$ for the macroscopic FPU chain:
(1) Replace Lax-shocks by dispersive shocks in the definition of
$\calS_\pm^\cons[u_\lst]$, i.e., use $\calD_\pm[u_\lst]$ from
\S\ref{s:RieSolDispShock}.
(2) Replace $r_2$ by the solution $\hat{r}_2$ of
$|c_\rmf(r_\lst,\hat{r}_2)| = |c_\rmrh(r_\lst,r_0)|$, and $r_1$ by
the solution $\hat{r}_1$ of $|c_\rmf(r_m,r_\lst)| =
|c_\rmrh(\hat{r}_1,r_m)|$.
\end{conjecture}
%
%
\subsection{Towards an FPU-Riemann solver for the non-classical
subsonic case} \label{s:slowSolver}
%
For subsonic conservative shock data, the situation is entirely
different. Recall that in Figure~\ref{f:consShock.B} subsonic
conservative shock initial data did not produce a conservative
shock. From \S\ref{s:tw} recall that fronts of the infinite chain do
not bifurcate in the subsonic case. Indeed, simulations analogous to
those in the supersonic case yield the drastically different results
plotted in Figures~\ref{fig:qp2_shock_curve} and
\ref{fig:qp2_rarefaction_curve}.  Since conservative shocks are
absent and composite wave of rarefaction and dispersive shocks
occur, we compare these with the solutions of the dissipative
p-system solver. It is the only solver in the (natural) family of
solvers studied in \cite{LeF02} that does not use non-classical
shocks, compare Figure~\ref{f:DissSolver.pSys}.
\bigpar
We investigate the modified wave sets for potential
\eqref{e:slowpot} as before. For given left state
$u_\lst=\pair{4}{0}$ we choose several values of $r_\rst$, see
Figure \ref{f:consdata}, and determine $v_\rst$ such that
$u_\rst=\pair{r_\rst}{v_\rst}$ belongs to $\calR_-[u_\lst]$ and
$\calS_+[u_\lst]$. The numerical results are plotted in Figures
\ref{fig:qp2_shock_curve} and \ref{fig:qp2_rarefaction_curve}; for
comparison with the dissipative p-system solver we use the inset
titles from Figure~\ref{f:DissSolver.pSys}.
%
%
For the shock initial data $S1$--$S5$ the chain produces single
dispersive shocks ($S1$--$S5\cong$ C) with increasing amplitudes,
decreasing back speeds $c_\rmb$ and increasing front speeds
$c_\rmf$. Here $c_\rmf$ is always larger than the corresponding
Rankine-Hugeniot speed $c_\rmrh$ and the speed of the conservative
shock. For $S6$--$S8\cong$ C+R we find a qualitatively different
solution with increasing rarefaction waves that are attached to the
same dispersive shock.
\begin{figure}[t!]%
\centering{%
\includegraphics[width=\mhpicDwidth, draft=\figdraft]%
{\figfile{qp2_shock_curve/qp2_sc_dist_mi_S1.eps}}%
\includegraphics[width=\mhpicDwidth, draft=\figdraft]%
{\figfile{qp2_shock_curve/qp2_sc_dist_mi_S3.eps}}%
\includegraphics[width=\mhpicDwidth, draft=\figdraft]%
{\figfile{qp2_shock_curve/qp2_sc_dist_mi_S5.eps}}%
\vskip\mhpicDvsep%
\includegraphics[width=\mhpicDwidth, draft=\figdraft]%
{\figfile{qp2_shock_curve/qp2_sc_dist_mi_S6.eps}}%
\includegraphics[width=\mhpicDwidth, draft=\figdraft]%
{\figfile{qp2_shock_curve/qp2_sc_dist_mi_S7.eps}}%
\includegraphics[width=\mhpicDwidth, draft=\figdraft]%
{\figfile{qp2_shock_curve/qp2_sc_dist_mi_S8.eps}}%
}%
\caption{%
  Simulations for data on the $2$-shock curve for potential
  \eqref{e:slowpot} and $u_\lst=(4,0)$, and $\MaLagr_\ast=0.1$,
  $N=3000$, $\MaTime=0.4$. The values for $r_\rst$ are those in
  Figure~\ref{f:consdata}(b), so `IVP $Sx$' crosses the turning point
  for increasing $x$.  The vertical lines indicate the locations of the
  $2$-shocks of the p-system corresponding to $\pair{r_\lst}{r_\rst}$,
  and the conservative shock corresponding to $u_\lst$.  }%
\label{fig:qp2_shock_curve}%
\quad\newline%
\centering{%
\includegraphics[width=\mhpicDwidth, draft=\figdraft]%
{\figfile{qp2_rarefaction_curve/qp2_rc_dist_mi_R2.eps}}%
\includegraphics[width=\mhpicDwidth, draft=\figdraft]%
{\figfile{qp2_rarefaction_curve/qp2_rc_dist_mi_R3.eps}}%
\includegraphics[width=\mhpicDwidth, draft=\figdraft]%
{\figfile{qp2_rarefaction_curve/qp2_rc_dist_mi_R4.eps}}%
\vskip\mhpicDvsep%
\includegraphics[width=\mhpicDwidth, draft=\figdraft]%
{\figfile{qp2_rarefaction_curve/qp2_rc_dist_mi_R5.eps}}%
\includegraphics[width=\mhpicDwidth, draft=\figdraft]%
{\figfile{qp2_rarefaction_curve/qp2_rc_dist_mi_R6.eps}}%
\includegraphics[width=\mhpicDwidth, draft=\figdraft]%
{\figfile{qp2_rarefaction_curve/qp2_rc_dist_mi_R7.eps}}%
}%
\caption{%
  Simulations for data on the $1$-rarefaction curve for
  potential \eqref{e:slowpot} analogous to Figure
  \ref{fig:qp2_shock_curve}.
}%
\label{fig:qp2_rarefaction_curve}%
\end{figure}
\par%
On the other hand, in the sequence $R1$--$R8$ the solutions are
single rarefaction waves ($R2\cong$ R), from which a dispersive
shock nucleates between $R3$ and $R4 \cong $ R+N, when crossing the
turning point $r_*$, see Figure \ref{f:consdata}(b). The rarefaction
fan shrinks from $R5$ to $R7$, and eventually the solution consists
of a single dispersive shock (not shown) corresponding to inset C in
Figure~\ref{f:DissSolver.pSys}(b). We conjecture that the solutions
from both Figures \ref{fig:qp2_shock_curve} and
\ref{fig:qp2_rarefaction_curve} can be understood by modifying the
dissipative solver as explained below. However, since there are no
non-classical shocks to test against, and since we cannot compute
$c_\rmf$ from the left and right states, our evidence is weaker than
in the supersonic case.
\bigpar%
\begin{figure}[t!]
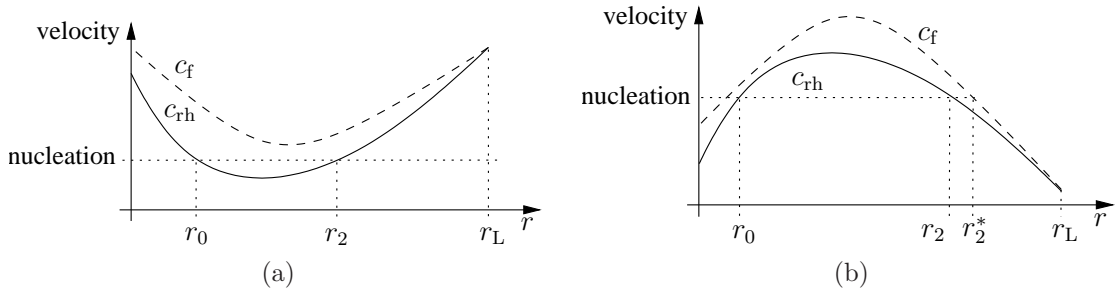

\centering{%
\begin{tabular}{cc}%
\includegraphics[width=0.45\textwidth, draft=\figdraft]%
{\figfile{xfig/cons_shock/slow-shock-vels}}%
&%
\includegraphics[width=0.45\textwidth, draft=\figdraft]%
{\figfile{xfig/cons_shock/fast-shock-vels}}%
\\(a) & (b)%
\end{tabular}%
}%
\caption{(a) Sketch of the explanation for non-nucleation
in
  the subsonic case. Here the front velocity $c_\rmf$ is larger than the
  nucleation velocity so that conservative shocks cannot be
  selected. (b) In the supersonic case, the nucleation cannot be missed in
  this way and occurs here at $r_2^*$.}
\label{f:non-nucleate}%
\end{figure}
We numerically observed the absence of subsonic conservative shocks
for various potentials. An explanation on the level of the Riemann
solver is the following.  Along the conservative shock-curve
$\calS_+^\cons[u_\lst]$ of the p-system, the nucleation of the
conservative shock connecting to $u_\lst$ occurs when it is as slow
as the compressive shock, i.e.,
$c_\rmrh(r_\lst,r_0)=c_\rmrh(r_\lst,r_\rst)$.  However, in the FPU
case, the criterion is naturally modified to equality of
conservative shock velocity and front velocity of the dispersive
shock, i.e., $c_\rmrh(r_\lst,r_0)=c_\rmf(r_\lst,r_\rst)$. Recall
that we have $c_\rmf(r_\lst,r_\rst) > c_\rmrh(r_\lst,r_\rst)$ for
$r_\rst\neq r_\lst$ due to Observation \ref{o:dispShock1}, and note
that both velocities converge to the characteristic velocity as
$r_\rst\to r_\lst$. It is therefore plausible that
$c_\rmf(r_\lst,r_\rst)>c_\rmrh(r_\lst,r_0)$ for all $r_\rst$ so that
the nucleation criterion always fails. We sketch this situation in
Figure~\ref{f:non-nucleate}(a). Note that the relative locations of
the $2$-shock and the dispersive shock front in
Figure~\ref{fig:qp2_shock_curve} support that this ordering indeed
occurs for the potential \eqref{e:slowpot}.
\par%
In contrast, in the supersonic case $\Phi^{(4)}(r_*)<0$ the velocity
curve is unimodal with a \emph{maximum}, so that nucleation cannot
be missed in this way, see Figure~\ref{f:non-nucleate}(b). We thus
arrive at the following conjecture.
\begin{conjecture}
  Let $r_*$ be the unique turning point of $\Phi'$. If
  $\Phi^{(4)}(r_*)<0$, then non-classical shocks are absent and the
  solution is qualitatively according to the dissipative solver. More
  precisely, the wave sets depicted in Figure~\ref{f:DissSolver.pSys} need
  to be modified as follows: (1) Replace Lax-shocks by dispersive shocks in the
  definition of $\calS_\pm^\cons[u_\lst]$, i.e., use $\calD_\pm[u_\lst]$
  from \S\ref{s:RieSolDispShock}.
  (2) Replace $r_0^*$ by the solution $\hat{r}_0^*$ to
  $c_\rmf(\hat{r}_0^*,r_\lst)^2=\Phi''(\hat{r}_0^*)$ and $r_1^*$ by
  the solution $\hat{r}_1^*$ to $c_\rmf(\hat{r}_0^*,r_\lst)^2=\Phi''(r_\lst)$.
\end{conjecture}
%
%
%
%
%
\section{Properties of conservative shocks}\label{s:mathconsshock}
\begin{figure}[t!]
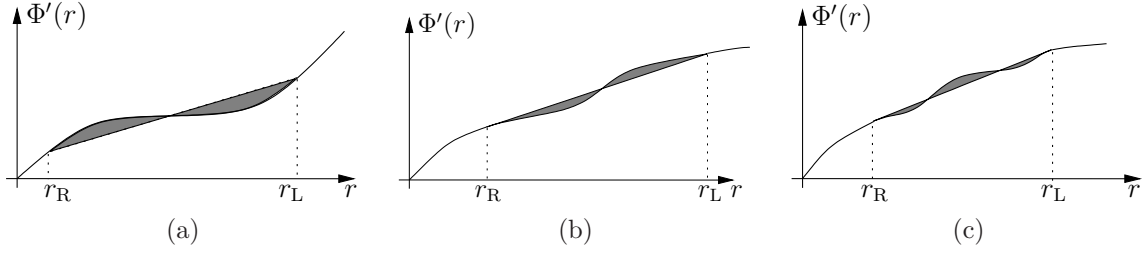

\centering%
\begin{tabular}{ccc}%
\includegraphics[width=0.3\textwidth, draft=\figdraft]%
{\figfile{xfig/cons_shock/consRel-ok}}%
&%
\includegraphics[width=0.3\textwidth, draft=\figdraft]%
{\figfile{xfig/cons_shock/consRel}}%
&%
\includegraphics[width=0.3\textwidth, draft=\figdraft]%
{\figfile{xfig/cons_shock/consRel-extr2}}%
\\(a) & (b) & (c)%
\end{tabular}%
\caption{Sample sketches for $\calJ(r_\lst,r_\rst)=0$. Between
  $r_\lst$ and $r_\rst$ the areas above and below the secant line up
  to the graph of $\Phi'$ (shaded) are equal. (a) The secant
  transversely intersects the graph on left and right so that
  the set $D$ is locally a monotone curve in the
  $(r_\lst,r_\rst)$-plane. (b) Secant and graph are tangent at
  $r_\lst$ so that $D$ has a local extremum in the
  $(r_\lst,r_\rst)$-plane. $\Phi'$ has (at least) two turning
  points in $[r_\rst,r_\lst]$. (c) The secant is tangent at both
  intersection points, and both tangencies point in the same
  direction, hence it is a local extremum of $\calJ$ and the point $(r_\lst,r_\rst)$ is isolated in $D$. $\Phi'$ has (at
  least) four turning points in $[r_\rst,r_\lst]$}
  \label{f:consGeom}
\end{figure}

%
In this section we study conservative shocks, that means solutions
to the three independent jump conditions \eqref{e:consJumpData}.
Eliminating the velocities $v_\lst$, $v_\rst$, an $c_\rmrh$ one
finds that each conservative shock is an element of
\begin{align}
\label{Eqn:ZeroJ} D=\left\{ \pair{r_\lst}{r_\rst}\;:\;
\calJ\pair{r_\lst}{r_\rst}=0\right\},
\end{align}
with $\calJ$ as in \eqref{e:ConsJump}. Conversely, each point in $D$
defines both a conservative $1$-shock and a $2$-shock, which are
unique up to Galilean transformation and differ only in
$\sgn\jump{v}=\sgn\,{c_\rmrh}$. The geometric interpretation of
$\calJ=0$ is that the signed area between the graphs of $\Phi'$ and
the secant line through $\Phi'(r_\lst)$ and $\Phi'(r_\rst)$ is zero
over $[r_\lst,r_\rst]$, compare Figure~\ref{f:consGeom}.
\bigpar

To characterise the structure of $D$ in presence of several
turning points of $\Phi^\prime$ we use the notation
$c_\lst:=|\lambda_\pm(r_\lst)|$ and $c_\rst:=|\lambda_\pm(r_\rst)|$.
\begin{theorem}
\label{t:consData}%
For $\Phi\in C^4(\R)$ the set $D$ has the following
properties.
\enum
\item
  Off-diagonal data $(r_\lst,r_\rst)\in D$ with
  $r_\lst\neq r_\rst$ exists if and only if $\Phi'$ has at least one
  turning point in the interval $(r_\lst,r_\rst)$.
\item
  Let $I\subset \R$ be any interval containing a
  \emph{single} turning point of $\Phi'$. Then $D\cap I\times{I}$ is the graph of a
  strictly decreasing function which crosses the diagonal $\{r_\lst=r_\rst\}$
  at the turning point.
\item
  The conservative shocks corresponding to $D\cap{I}\times{I}$ are
  undercompressive. If $\Phi^{(4)}(r_*)<0$, they are supersonic, and
  if $\Phi^{(4)}(r_*)>0$ subsonic.
\item\label{i:compr}
  Compression changes precisely at local extrema
  in the coordinate directions of $D$ in the
  $(r_\lst,r_\rst)$-plane. At extrema in the $r_\lst$-direction
  $c_\rst^2$ crosses $c_\rmrh(r_\lst,r_\rst)^2$, and at extrema in the
  $r_\rst$-direction $c_\lst^2$ crosses $c_\rmrh(r_\lst,r_\rst)^2$.
\item\label{i:conn}
  If $\Phi'$ has precisely two turning points, then
  $D$ is the union of the diagonal $\{r_\lst=r_\rst\}$ and a closed
  curve crossing the diagonal at the turning points.
\item
  The set $D$ does not have bounded connected components
  if $\Phi'$ has three or fewer turning points, and
  $D\setminus\{r_\lst=r_\rst\}$ is bounded if the number of turning
  points is even.
\end{list}
\end{theorem}
\begin{figure}
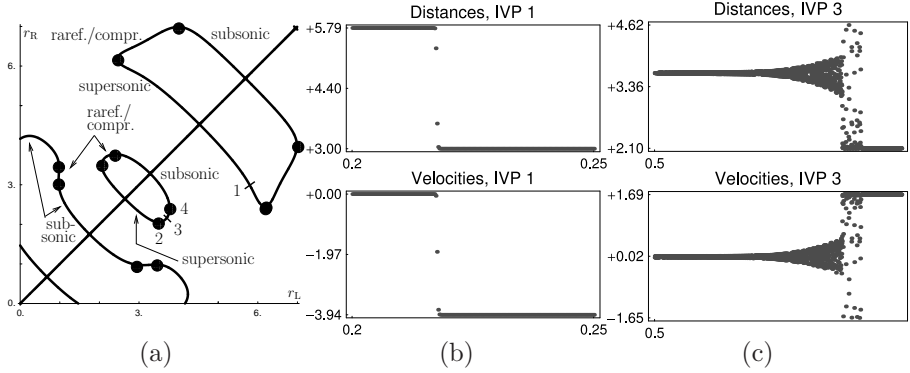

\centering{ %
\setlength{\tabcolsep}{0cm}%
\begin{tabular}{ccc}%
\begin{minipage}[c]{\mhpicDwidth}%
\includegraphics[width=\mhpicDwidth, draft=\figdraft]%
{\figfile{xfig/ivp_data/pot_sc_ivp_data.eps}}%
\end{minipage}%
&%
\begin{minipage}[c]{\mhpicDwidth}%
\includegraphics[width=\mhpicDwidth, draft=\figdraft]%
{\figfile{pot_sc/sc_pnt_1_dist}}%
\\%
\includegraphics[width=\mhpicDwidth, draft=\figdraft]%
{\figfile{pot_sc/sc_pnt_1_vel}}%
\end{minipage}%
&%
\begin{minipage}[c]{\mhpicDwidth}%
\includegraphics[width=\mhpicDwidth, draft=\figdraft]%
{\figfile{pot_sc/sc_pnt_5_dist}}%
\\%
\includegraphics[width=\mhpicDwidth, draft=\figdraft]%
{\figfile{pot_sc/sc_pnt_5_vel}}%
\end{minipage}%
\\(a) & (b) & (c)%
\end{tabular}%
}%
\caption{%
  (a) The set $D$ for the potential \eqref{e:potSC}, i.e., on solid
  curves holds $\calJ=0$; off-diagonal segments of $D$ are labelled according to
  the ordering of $c_\rmrh^2$ and $\lambda_\pm^2$. (b) Solution for a
  conservative shock with initial data at point (1) in (a). (c)
  Solution for the initial data from point (3) in (a); the
  solutions appears to contain a jump from wave train to wave
  train. In (b) \& (c) the initial data ranges over two turning points
  of $\Phi'$; $N=4000$ initial jump $\bal_*=0.5$.} %
\label{f:potSC-cons-curves}%
\end{figure}%
In Figure~\ref{f:potSC-cons-curves}(a) part of the set $D$ is plotted
for the potential
\begin{equation}\label{e:potSC}
\Phi(r+1) = r  + \frac 1 2r^2 + \frac 1 {20}r^3 - \frac 1 4 \cos(2r)
+ \frac 1 {10} \sin(3r ).
\end{equation}
Changes in shock type occur for instance at point 2, which is an
extremum in the $r_\rst$-direction so that $c_\lst^2$ becomes larger
than $c_\rmrh^2$ in the direction towards the nearest extremum in
$r_\lst$-direction at point 4.  The conservative shocks on this curve
start out supersonic, hence in the segment between points 2 and 4 the
1-shocks are compressive and the 2-shocks are rarefaction shocks.  At
point 4 the term $c_\rst^2$ becomes larger than $c_\rmrh^2$ when
crossing it away from point 3, and so conservative shocks beyond point
4 are subsonic.
\begin{remark}
\label{Rem:Prop.Cons.Shocks.1}
\enum
\item
  Conservative shocks do not have a preferred direction of
  propagation, and are isolated points on the shock curves
  $\calS_\pm[u_\lst]$. In contrast, classical shocks have a selected direction of propagation in
  order to be compressive and generate continuous segments of
  $\calS_\pm[u_\lst]$.
\item%
  In lack of turning points, conservative shock data does not exist
  for the Toda potential and any cubic
  potential. For harmonic potentials, however, all shocks are
  conservative, and in fact contact discontinuities.
\item%
  For even potentials off-diagonal conservative jumps occur for
  $r_\lst=-r_\rst$, because $\calJ(r,-r)\equiv 0$ by symmetry. More generally,
  the symmetry $\Phi(r_*+r)=\Phi(r_*-r)$ with $\Phi'''(r_*)=0$ provides
  $\tilde{D}\subset{D}$ with
  \begin{math} %
  \tilde{D}=\left\{(r_\lst,r_\rst)=(r_*+r, r_*-r):r\in\R\right\}.
  \end{math} %
  Note that each conservative shock from $\tilde{D}$ is degenerate as
  the characteristic velocities for left and
  right states equal.
\item%
  For quartic potentials (i.e., the classical FPU chains) all
  off-diagonal conservative data are given by $\tilde{D}$, so that
  non-degenerate conservative shocks occur only for potentials of
  polynomial degree five or higher.
\end{list}
\end{remark}
\bigpar
\begin{Proof}[~of Theorem~\ref{t:consData}]
\enum
\item
  This immediately follows from the geometric interpretation of
  $\calJ=0$, see Figure~\ref{f:consGeom}.
\item Again, the geometric interpretation shows that for fixed
  $r_\lst\subset I$ there is at most one solution $r\in I$ to
  $\calJ(r_\lst,r)=0$; similar for $r_\rst$. Hence, the solution set
  is a monotone curve; that it decays follows similarly. Note that
  tangents of $\Phi$ and the secant slope cannot coincide within $I$.
\item \addtocounter{Lcount}{1} \& \arabic{Lcount}. 
These follow from Lemma~\ref{l:consShock} below.
%
%
\item
  Consider $(r_\lst,r_\rst)$ as in Figure \ref{f:consGeom}(b), where
  $\calJ(r_\lst,r_\rst)=0$ and the secant slope coincides with the
  tangent slope of $\Phi'$ at one end. Note that, since there are only
  two turning points, the graph of $\Phi'$ must lie on one side of the
  secant line near the tangency. Moving monotonically through the
  point where these slopes coincide, the area between the graphs can
  only vanish when the other point reverses its direction.  Hence, the
  curve in the $(r_\lst,r_\rst)$-plane has an extremum, and this can
  only occur when both turning points lie in the interval
  $(r_\rst,r_\lst)$.  When continuing the curves from item 2, the
  tangency points of the graph intersections must be reached. Since
  there are no further changes in monotonicity, and curves are unique
  in suitable intervals $I$, the two curves emanating from the turning
  points must connect.

  It remains to show that there can be no bounded components of $D$
  that are isolated from the diagonal.

Note that a stationary point $\pair{r_\lst^*}{r_\rst^*}$ of $\calJ$
requires tangency of the graph of $\Phi'$ with the secant segment at
\emph{both} $r_\lst^*$ and $r_\rst^*$.  Such a stationary point is a
local extremum if the graph of $\Phi'$ is either above or below the
common tangency at both $r_\lst^*$ and $r_\rst^*$, i.e.,
$\Phi'''\at{r_\lst^*}$ and $\Phi'''\at{r_\rst^*}$ have the same
sign. For $\Phi$ with two turning points there exist a unique
stationary point which is moreover a local extremum
$(r^*_\lst,r^*_\rst)$ with $\calJ(r^*_\lst,r^*_\rst)\neq 0$, because the
geometry implies that the enclosed area is only on one side of the
secant.  Hence, each zero of $\calJ$ is a regular point and cannot be
an isolated point of $D$. Now suppose for contradiction that a
connected and bounded component of $D$ existed.  Then it must be a
closed curve containing the local extremum $\pair{r_\lst^*}{r_\rst^*}$
in its interior.  However, fixing $r^*_\lst$ and moving from
$r^*_\rst$ towards $r^*_\lst$ the secant segment stays above or below
the graph of $\Phi'$ so that $\calJ\neq 0$ until $r_\lst=r_\rst$ on
the diagonal, which is the contradiction.
\item
    We continue the discussion of bounded isolated components from the
    previous item with a local extremum $\pair{r_\lst^*}{r_\rst^*}$ in
    the interior.  Our arguments from above imply that the interval
    $[r_\rst^*,r_\lst^*]$ contains at least three turning points, and
    the tangency criterion for local extrema shows that the number of
    turning points in the interval must be even.

  Concerning boundedness of $D\setminus\{r_\lst=r_\rst\}$, note that
  for an even number of turning points, the convexity of $\Phi'$
  outside a sufficiently large interval is the same. Hence, the secant
  line for sufficiently far distant $r_\lst$ and $r_\rst$ lies on one
  side of the graph of $\Phi'$ so that $\calJ\neq 0$.
\end{list}
\end{Proof}
Using the proof of the last item, it is not difficult to construct
$\Phi$ for which the set $D$ consists of several bounded components
that are disconnected from each other. The following lemma gives
some more specific information.
\begin{lemma}\label{l:consShock}~
\enum
\item Whenever $\Phi'''(r_*)=0$ and $\Phi^{(4)}(r_*)\neq 0$ for some
  $r_*\in \R$ there exist a smooth locally unique curve
  $r_\lst\mapsto{R}\at{r_\lst}$ of solutions to $\calJ=0$ in
  $\{r_\lst\neq r_\rst\}\cup\{(r_*,r_*)\}$ and it has tangent $(-1,1)$
  at $(r_*,r_*)$.
\item Whenever $c_\rmrh^2(r_\lst,r_\rst) \neq c_\rst^2$ and
  $\calJ(r_\lst,r_\rst)=0$ then the set $D$ from
  Theorem~\ref{t:consData} is locally given by a function
  $r_\rst=R(r_\lst)$, which has a local extremum if
  $c_\rmrh^2(r_\lst,r_\rst) = c_\lst^2$. Similarly, whenever
  $c_\rmrh^2(r_\lst,r_\rst) \neq c_\lst^2$ and
  $\calJ(r_\lst,r_\rst)=0$ then $D$ is locally given by a function
  $r_\lst=R(r_\rst)$, which has a local extremum if
  $c_\rmrh^2(r_\lst,r_\rst) = c_\rst^2$.
\item
On a curve $(r,R(r))$ as in (1.) it holds that
\begin{eqnarray*}
\sgn(c_\rmrh^2 - c_\pm^2) = -\sgn\nat{\Phi^{(4)}(r_*)}%
\text{ for  }r\approx r_*,\qquad \sgn(c_\lst^2 - c_\rst^2) =
+\sgn\nat{\Phi^{(5)}(r_*)}\text{ for }r\approx r_*.
\end{eqnarray*}
\end{list}
\end{lemma}
\begin{Proof}
\enum
\item
We readily compute that all first and second order partial
derivatives of $\calJ$ vanish on the diagonal $\{r_\lst=r_\rst\}$
and that all third order derivatives contain the factor
$\Phi'''(r_{\lst/\rst})$. Implicit differentiation of
$\calJ(r,R(r))$ with respect to $r$ then shows that bifurcations of
solutions to $\calJ=0$ from the diagonal can only occur for
$\Phi'''(r)=0$. The resulting bifurcation equation at $r_*$, using
$x=R'(r_*)$, reads
\begin{align*}
\Phi^{(4)}(r_*)(2x^4 -3x^3+3x-2)=0
\end{align*}
and has the solution $x=1$, corresponding to the trivial solution
curve along the diagonal, and $x=-1$, corresponding to the new
bifurcating branch, as well as two complex conjugate roots that do
not contribute to real solutions.
\item
Labelling variables as $\calJ(r_1,r_2)=\calJ(r_\lst,r_\rst)$ we
compute
\begin{align*}
\partial_{r_j}\calJ(r_1,r_2) = \frac{\Phi'(r_1)-\Phi'(r_2)}{2} - \frac{r_1-r_2}{2}\Phi''(r_j).
\end{align*}
Using this and the definitions of the velocities, implicit
differentiation of $\calJ(r,R(r))=0$ gives
\begin{align*}
R'(r) = -
\frac{\partial_{r_1}\calJ(r,R(r))}{\partial_{r_2}\calJ(r,R(r))}
=-\frac{c_\rmrh(r_\lst,r_\rst)^2-c_\lst^2}{c_\rmrh(r_\lst,r_\rst)^2-c_\rst^2}.
\end{align*}
The statement immediately follows from this formula.
\item
It follows from (1.) that $R(r_*+s)= r_*-s$ to leading order, so
that we can expand $G_\pm(s):=c_\rmrh^2(s_*+s,r_*-s) -
c_\pm^2(r_*\mp s)$ and to leading order we obtain
\begin{eqnarray*}
G_\pm(s) &=& \frac{\Phi'(r_*+s) - \Phi'(r_*-s)}{r_*+s - (r_*-s)} - \Phi''(r_* \mp s)\\
&=& \frac{2\Phi''(r_*)s + 2\Phi^{(4)}(r_*)\frac{s^3}{6}}{2s} -
\Phi''(r_*) - \Phi^{(4)}(r_*)\frac{s^2}{2} = -\frac{1}{3}s^2
\Phi^{(4)}(r_*)
\end{eqnarray*}
Therefore $\sgn(c_\rmrh^2 - c_\pm^2) = -\sgn(\Phi^{(4)}(r_*))$ for
$s\approx 0$. The proof for $c_\lst^2-c_\rst^2$ is completely analogous.
\end{list}
\end{Proof}
The results of this section can be readily generalised to degenerate
situations where $\Phi^{(4)}= 0$ at the turning
point: Theorem~\ref{t:consData} is unchanged, but the bifurcation
equations and local set of solutions in Lemma~\ref{l:consShock} does
according to the degree of degeneracy.
\par
Finally, concerning FPU solutions for more than one turning point,
consider the results for initial conservative 2-shock data that
ranges over two turning points in
Figure~\ref{f:potSC-cons-curves}(b),(c). The solution in (b) is the
expected conservative shock, but the one in (c) is not. In fact, the
solution in (c) appears to contain a (conservative) shock between
wave trains which cannot be predicted from the p-system at all.
\par
Thus, the predictive power of the p-system investigated in this
paper ends even qualitatively for more complicated fluxes. Extended
systems such as the Whitham modulation equations are required to
understand the solution structures, but as even the hyperbolic
nature of these is unknown, it is left for the future.
\appendix
%
%
\section{Classical Riemann solver for the p-system}\label{app:p-riemann}
%
%
Formal substitution into \eqref{e:ColdModEqn} of a self-similar ansatz
in the variable $c=\at{\bal-\bal_*}/\bt$ gives
\begin{equation}%
\label{e:pSys-c} -c \dot r = \dot v\;,\quad -c\dot v = \Phi''(r)\dot
r,
\end{equation}
where $\dot~=\rmd/\rmd c$, and this implies $c^2 =\Phi''(r)$. The
eigenvalues, i.e., characteristic velocities, and associated right
eigenvectors of the p-system are given by $\lambda_\pm(r) =
\pm\sqrt{\Phi''(r)}$ and $e_\pm(r) = (1,\lambda_\mp)$, respectively.
Strict convexity of $\Phi$ implies $\lambda_-(r)<0<\lambda_+(r)$
for all $r$, and the p-system is therefore (globally) strictly
hyperbolic. Eigenvalues are genuinely nonlinear as long as $\Phi'''$
does not change sign.
\par
Concerning symmetries, the p-system respects Galilean transformations,
and so the set of self-similar solutions is invariant under
\begin{math}
\pair{r(c)}{v(c)}\mapsto\pair{r(c)}{v(c)+v_0}
\end{math} %
where $v_0$ is constant. Moreover, the p-system exhibits the
\emph{reflection symmetry} that under $c\mapsto{-c}$ each
self-similar solution transforms according to
\begin{math}
\bpair{r(c)}{v(c)}\mapsto\bpair{r(-c)}{-v(-c)}.
\end{math} %
\par
The Lax theory for strictly hyperbolic systems with genuinely
nonlinear eigenvalues is built up from the following two types of
elementary waves. \emph{Rarefaction fans} are smooth solutions, and
\eqref{e:pSys-c} implies the existence of two families, called
$1$-waves and $2$-waves with $c=\lambda_-\at{r\at{c}}$ and
$c=\lambda_+\at{r\at{c}}$, respectively. \emph{Shocks} are
discontinuous solutions, and satisfy \eqref{e:pSys-c} only in the
sense of distributions. This gives rise to the Rankine-Hugeniot
conditions \eqref{e:ColdDataJump} with \emph{shock speed} $c_\rmrh$.
The convexity of $\Phi$ and the identity
\begin{math}
c_\rmrh^2\jump{r} = \jump{\Phi'(r)}
\end{math}
imply the existence of $1$-shocks with $c_\rmrh<0$ and $2$-shocks with
$c_\rmrh>0$. Lax theory considers only \emph{compressive} shocks that
satisfy the \emph{Lax condition} (with $-$ and $+$ for $1$- and
$2$-shocks, respectively)
\begin{align}
\label{eqn:LaxCond}
\lambda_\pm(u_{\lst})>c_\rmrh>\lambda_\pm(u_{\rst}),
\end{align}
\par
To describe the wave set $\calW[u_\lst]$ for given left state
$u_{\lst}=\pair{r_{\lst}}{v_{\lst}}$, we define the \emph{integral
  curves} $\calO_-[u_{\lst}]$ and $\calO_+[u_{\lst}]$ by
\begin{align*}
\pair{r}{v}\in\calO_\pm[u_{\lst}]\quad\Leftrightarrow\quad
v=v_{\lst} \mp \int_{r_{\lst}}^{r} \sqrt{\Phi''(s)}\rmd s,
\end{align*}
and the \emph{Hugeniot curves} $\calH_-[u_{\lst}]$ and
$\calH_+[u_{\lst}]$ given by
\begin{align*}
\pair{r}{v}\in\calH_\pm[u_{\lst}]\quad\Leftrightarrow\quad
v=v_{\lst}
\mp\sqrt{\at{\Phi'\at{r}-\Phi'\at{r_{\lst}}}\at{r-r_{\lst}}}.
\end{align*}
All these curves contain the point $u_{\lst}$ and can be globally
parameterized by $r$. Due to \eqref{e:pSys-c} we find
$\tfrac{\dint{v}}{\dint{r}}=-c$ which implies 
\begin{equation}
\notag
\pair{r_{\rst}}{v_{\rst}}\in\calR_\pm[u_{\lst}]\quad\Leftrightarrow\quad
\pair{r_{\rst}}{v_{\rst}}\in\calO_\pm[u_{\lst}]\quad\text{and}\quad
{\lambda_\pm\at{r_{\rst}}>\lambda_\pm\at{r_{\lst}}}.
\end{equation}
\par%
The \emph{shock wave set}
$\calS[u_{\lst}]=\calS_-[u_{\lst}]\cup\calS_+[u_{\lst}]$ consists of
all possible right states $u_\rst$ that can be connected with
$u_{\lst}$ by a single Lax shock. By construction, any $u_\rst$ must
be an element of one of the Hugeniot curves of $u_\lst$, but the
Lax-condition \eqref{eqn:LaxCond} selects one branch for each
Hugeniot curve.
\par
Exchanging the role of left and right states provides sets
$\widetilde{\calR}[u_{\rst}]=\widetilde{\calR}_-[u_{\rst}]\cup\widetilde{\calR}_+[u_{\rst}]$
and
$\widetilde{\calS}[u_{\rst}]=\widetilde{\calS}_-[u_{\rst}]\cup\widetilde{\calS}_+[u_{\rst}]$,
which contain all possible left states that can be connected to a
prescribed right state $u_{\rst}$ by a single rarefaction wave or
Lax shock, respectively. The standard Riemann solver is defined by
these left and right wave sets as follows, see Figure
\ref{fig:atomistic_rs}. For given $\pair{u_{\lst}}{u_{\rst}}$ there
exists a unique (see, e.g., \cite{Smo94}) intermediate state
$u_{\mst}\in\calW[u_{\lst}]\cap\widetilde{W}[u_{\rst}]$ such that
$u_{\mst}\in{\calW_-[u_{\lst}]}$ and
$u_{\rst}\in{\calW_+[u_{\mst}]}$. The solution of the Riemann
problem consists of the two elementary waves that connect $u_\lst$
to $u_\mst$, and $u_\mst$ to $u_\rst$ (one of these may be trivial).
\par
Shocks are classified as follows. A shock connecting $u_\lst$ to
$u_\rst$ with speed $c_\rmrh$ is called
\begin{enumerate}
\item
\emph{compressive}, or \emph{Lax shock}, if
$\lambda\at{u_{\lst}}>c_\rmrh>\lambda\at{u_{\rst}}$,
\item
\emph{rarefaction shock}, if
$\lambda\at{u_{\lst}}<c_\rmrh<\lambda\at{u_{\rst}}$,
\item
\emph{supersonic}, or \emph{fast undercompressive}, if
$\abs{c_\rmrh}>\abs{\lambda\at{u_{\lst}}}$ and
$\abs{c_\rmrh}>\abs{\lambda\at{u_{\rst}}}$,
\item
\emph{subsonic}, or \emph{slow undercompressive}, if
$\abs{c_\rmrh}<\abs{\lambda\at{u_{\lst}}}$ and
$\abs{c_\rmrh}<\abs{\lambda\at{u_{\rst}}}$,
\item
\emph{sonic}, if $\abs{c_\rmrh}=\abs{\lambda\at{u_{\lst}}}$ or
$\abs{c_\rmrh}=\abs{\lambda\at{u_{\rst}}}$,
\end{enumerate}
with $\sgn\,{c}_\rmrh=\sgn\,{\lambda}<0$ and
$\sgn\,{c}_\rmrh=\sgn\,{\lambda}>0$ for $1$- and $2$- shocks,
respectively. All these definitions are invariant under reflections
$c_\rmrh\leftrightsquigarrow-c_\rmrh$,
$\lambda\at{u_\rst}\leftrightsquigarrow-\lambda\at{u_\lst}$.
\section*{Acknowlegdements}%
This work has been supported in part by the DFG Priority Program
1095 ``Analysis, Modeling and Simulation of Multiscale Problems''
(M.H., J.R.), and the NDNS+ cluster of the NWO (J.R.). We thank
Wolfgang Dreyer for motivating and supporting part of this research,
and Alexander Mielke as well as Thomas Kriecherbauer for fruitful
discussions. Finally, we are grateful to the editor for handling the 
paper professionally and the anonymous reviewers for
helping us to improve the exposition of the material.
\bibliographystyle{amsalpha}%
\providecommand{\bysame}{\leavevmode\hbox to3em{\hrulefill}\thinspace}
\providecommand{\MR}{\relax\ifhmode\unskip\space\fi MR }
\providecommand{\MRhref}[2]{%
  \href{http://www.ams.org/mathscinet-getitem?mr=#1}{#2}
}
\providecommand{\href}[2]{#2}

\end{document}